\documentclass[11pt,a4paper,reqno]{amsart}
\usepackage{amsfonts,amsthm, amscd, epsfig, amsmath, amssymb,enumerate}
\usepackage{amstext}
\usepackage{xspace}
\usepackage{url}
\usepackage{enumerate}
\usepackage[all,2cell,ps]{xy}
\usepackage{color}

\newcommand{\halmos}{\rule{1ex}{1.4ex}}

\newtheorem{itlemma}{Lemma}[section]
\newtheorem{itproposition}[itlemma]{Proposition}
\newtheorem{itfact}[itlemma]{Fact}
\newtheorem{theorem}[itlemma]{Theorem}
\newtheorem{itcorollary}[itlemma]{Corollary}
\newtheorem{itremark}[itlemma]{Remark}
\newtheorem{itremarks}[itlemma]{Remarks}
\newtheorem{itdefinition}[itlemma]{Definition}
\newtheorem{itexample}[itlemma]{Example}

\newenvironment{fact}{\begin{itfact}\rm}{\end{itfact}}
\newenvironment{claim}{\begin{itclaim}\rm}{\end{itclaim}}
\newenvironment{lemma}{\begin{itlemma}}{\end{itlemma}}
\newenvironment{remark}{\begin{itremark}\rm}{\end{itremark}}
\newenvironment{remarks}{\begin{itremarks} \rm}{\end{itremarks}}
\newenvironment{corollary}{\begin{itcorollary}}{\end{itcorollary}}
\newenvironment{proposition}{\begin{itproposition}}{\end{itproposition}}
\newenvironment{definition}{\begin{itdefinition}\rm}{\end{itdefinition}}
\newenvironment{example}{\begin{itexample}\rm}{\end{itexample}}
\newcommand{\be}[1]{\begin{equation}\label{#1}}
\newcommand{\ee}{\end{equation}}
\newcommand{\bl}[1]{\begin{lemma}\label{#1}}
\newcommand{\br}[1]{\begin{remark}\label{#1}}
\newcommand{\brs}[1]{\begin{remarks}\label{#1}}
\newcommand{\bt}[1]{\begin{theorem}\label{#1}}
\newcommand{\bd}[1]{\begin{definition}\label{#1}}
\newcommand{\bp}[1]{\begin{proposition}\label{#1}}
\newcommand{\bfact}[1]{\begin{fact}\label{#1}}
\newcommand{\bc}[1]{\begin{corollary}\label{#1}}
\newcommand{\bex}[1]{\begin{example}\label{#1}}
\newcommand{\ec}{\end{corollary}}
\newcommand{\efact}{\end{fact}}
\newcommand{\eex}{\end{example}}
\newcommand{\el}{\end{lemma}}
\newcommand{\er}{\end{remark}}
\newcommand{\ers}{\end{remarks}}
\newcommand{\et}{\end{theorem}}
\newcommand{\ed}{\end{definition}}
\newcommand{\ep}{\end{proposition}}
\newcommand{\epr}{\end{proof}}
\newcommand{\bpr}{\begin{proof}}
\newcommand{\bcl}[1]{\begin{claim}\label{#1}}
\newcommand{\ecl}{\end{claim}}

\newcommand{\ecs}{\end{corollary}}
\newcommand{\eers}{\end{exercise}}
\newcommand{\eexs}{\end{example}}
\newcommand{\eems}{\end{example}}
\newcommand{\els}{\end{lemma}}
\newcommand{\eles}{\end{lemmaex}}
\newcommand{\ets}{\end{theorem}}
\newcommand{\eds}{\end{definition}}
\newcommand{\eps}{\end{proposition}}

\newcommand{\bi}{\begin{itemize}}
\newcommand{\ei}{\end{itemize}}
\newcommand{\ben}{\begin{enumerate}}
\newcommand{\een}{\end{enumerate}}

\def\vbar{\mathchoice{\vrule height6.3ptdepth-.5ptwidth.8pt\kern-.8pt}
   {\vrule height6.3ptdepth-.5ptwidth.8pt\kern-.8pt}
   {\vrule height4.1ptdepth-.35ptwidth.6pt\kern-.6pt}
   {\vrule height3.1ptdepth-.25ptwidth.5pt\kern-.5pt}}
\def\fudge{\mathchoice{}{}{\mkern.5mu}{\mkern.8mu}}
\def\bbc#1#2{{\rm \mkern#2mu\vbar\mkern-#2mu#1}}
\def\bbb#1{{\rm I\mkern-3.5mu #1}}
\def\bba#1#2{{\rm #1\mkern-#2mu\fudge #1}}
\def\bb#1{{\count4=`#1 \advance\count4by-64 \ifcase\count4\or\bba A{11.5}\or
   \bbb B\or\bbc C{5}\or\bbb D\or\bbb E\or\bbb F \or\bbc G{5}\or\bbb H\or
   \bbb I\or\bbc J{3}\or\bbb K\or\bbb L \or\bbb M\or\bbb N\or\bbc O{5} \or
   \bbb P\or\bbc Q{5}\or\bbb R\or\bbc S{4.2}\or\bba T{10.5}\or\bbc U{5}\or
   \bba V{12}\or\bba W{16.5}\or\bba X{11}\or\bba Y{11.7}\or\bba Z{7.5}\fi}}

\def \qed {{\hspace*{\fill}$\halmos$\medskip}}

\newcommand {\cal}[1] {\mathcal{#1}}

\def \R {{\mathbb R}}

\def \N {{\mathbb N}}

\def \ra {\rightarrow }

\def \s{\sigma}

\def \P{{\mathbb{P}}}
\def\g{\gamma}
\def\d{\delta}
\def\e{\varepsilon}

\def\b{\beta}

\def\l{\lambda}
\def\E{{\mathrm{E}}}
\def\var{{\mathrm{Var}}}
\def\F{{\mathcal{F}}}

\def\1{{\bf 1}}
\def\e{\varepsilon}

\def\ec{\`e }

\begin{document}

\title[Synchro and FCLT for interacting RW]{Synchronization and
  Functional Central Limit Theorems \\for
Interacting Reinforced Random Walks}
  \maketitle

    \centerline{Irene Crimaldi\footnote{IMT School for Advanced
        Studies Lucca, Piazza San Ponziano 6, 55100 Lucca, Italy,
        \url{irene.crimaldi@imtlucca.it}} , Paolo Dai
      Pra\footnote{Dipartimento di Matematica, Universit\`a degli
        Studi di Padova, Via Trieste 63, 35121 Padova, Italy,
        \url{daipra@math.unipd.it}},
      Pierre-Yves~Louis\footnote{Laboratoire de Math\'ematiques et
        Applications UMR 7348, Universit\'e de Poitiers et CNRS,
        T\'el\'eport 2, BP 30179 Bvd Marie et Pierre Curie, 86962
        Futuroscope Chasseneuil Cedex, France,
        \url{pierre-yves.louis@math.cnrs.fr}}, Ida
      G.~Minelli\footnote{Dipartimento di Ingegneria e Scienze
        dell'Informazione e Matematica, Universit\`a degli Studi
        dell'Aquila, Via Vetoio (Coppito 1), 67100 L'Aquila, Italy,
        \url{ida.minelli@dm.univaq.it}}}

\begin{abstract}
We obtain Central Limit Theorems in Functional form for a class of
time-inhomogeneous interacting random walks on the simplex of
probability measures over a finite set. Due to a reinforcement
mechanism, the increments of the walks are correlated, forcing their
convergence to the same, possibly random, limit.  Random walks of this
form have been introduced in the context of urn models and in
stochastic algorithms. We also propose an application to opinion
dynamics in a random network evolving via preferential attachment.  We
study, in particular, random walks interacting through a mean-field rule
and compare the rate they converge to their limit with the rate of
{\em synchronization}, i.e. the rate at which their mutual distances
converge to zero. Under certain conditions, synchronization is faster
than convergence.
\end{abstract}



\section{Introduction}

Let $S$ be a finite set and denote by $\mathcal{P}(S)$ the simplex of
probabilities on $S$:
\[
\mathcal{P}(S) := \left\{ \mu:S \ra [0,1]: \ \sum_{x \in S} \mu(x) = 1 \right\}.
\]
In this paper we consider stochastic evolutions on $\mathcal{P}(S)$ of
the form \be{dyn-symplex} \mathcal{Z}_{n+1} = (1-r_n) \mathcal{Z}_n +
r_n K_n(I_{n+1}), \ee where $0 \leq r_n < 1$ are given numbers,
$K_n:S \ra \mathcal{P}(S)$ are given functions, and $(I_n)_{n \geq 1}$
is a sequence of $S$-valued random variables such that, for $\F_n :=
\s(\mathcal{Z}_0,\mathcal{Z}_1,\ldots,\mathcal{Z}_n)$, \be{rdmI}
\P(I_{n+1} = x | \F_n) = \mathcal{Z}_n(x).  
\ee 
We think of this as a \emph{generalized reinforcement mechanism}: note
indeed that, in the particular case $K_n(x)=\delta_x$, with $\delta_x$
denoting the Dirac measure at $x\in S$, the larger ${\mathcal
  Z}_n(x)$, the higher the probability of increasing it at the next
step.\\

\indent Models of type \eqref{dyn-symplex} can be viewed as
time-inhomogeneous random walks on $\mathcal{P}(S)$, and arise
naturally in at least two distinct contexts.  \\

\medskip

\noindent{\em 1.1. Urn Models.} Let $S$ be the set of the colors of
the balls in a urn. Consider the following scheme. A ball is randomly
drawn, uniformly among all balls. Suppose this is the $(n+1)$-st
draw. If its color is $y$ then we reinsert it in the urn and, for each
color $x \in S$, we add $A_n(y,x)$ balls of color $x$, where
\[
A_n: S \times S \ra \N_0 = \{0,1,2,\ldots\}
\]
is a given function, called {\em reinforcement matrix}, satisfying the
{\em balance} condition: the sum \\ ${\displaystyle{\sum_{x \in S}
    A_n(y,x) = \overline{A}_n}}$ does not depend on $y$. As a
consequence, the total number $$N(n) = N(0) + \sum_{k=0}^{n-1}
\overline{A}_k$$ of balls in the urn after $n$ steps does not depend
on the sequence of colors drawn. Denote by $N(n,x)$ the number of
balls in the urn of color $x$ after $n$ steps, and
\[
\mathcal{Z}_n(x) := \frac{N(n,x)}{N(n)}.
\]
Then $\mathcal{Z}_n \in \mathcal{P}(S)$, and it evolves as in
\eqref{dyn-symplex} with
\[
r_n :=
\frac{\overline{A}_{n}}{N(n+1)} \ \ \ K_n(y)(\cdot)
= \frac{A_n(y,\cdot)}{\overline{A}_n}.
\]
This model includes the P\'olya and the Friedman scheme as special
cases, as well as many generalizations with time dependent
reinforcement scheme (see e.g. \cite{mah} for an introduction to the
subject). Note that in the most classical schemes (P\'olya, Friedman)
$\overline{A}_n$ is constant in $n$. More generally, in all cases in
which $\overline{A}_n$ grows at most polynomially in $n$, we have that
$r_n$ is of order $\frac{1}{n}$ as $n \ra +\infty$.  \\

\medskip
\noindent{\em 1.2. Opinion dynamics on preferential attachment
  graphs.}  Consider a sequence of random non-oriented graphs $G_n =
(V_n, E_n)$, evolving through a preferential attachment rule (see
e.g. \cite{bar, remco}). More specifically, for a given $\d > -1$, the
graph evolves according to the following rules: \bi
\item
at time $n=2$ the graph consists of the two vertices $\{1,2\}$
connected by one edge;
\item
at time $n+1$ the new vertex $n+1$ is added and it is linked with an
edge to vertex $i \in V_n = \{1,2,\ldots,n\}$ with probability
$\frac{d_i(n) + \d}{2(n-1)+n\d}$, where $d_i(n)$ is the {\em degree}
of the vertex $i$ at time $n$, i.e. the number of edges having $i$ as
endpoint.  \ei Note that $G_n$ is a connected graph.

We now define a stochastic dynamics, whose evolution depends on the
realization of the graph sequence $(G_n)_{n \geq 2}$, which therefore
plays a role analogous of that of a dynamic {\em random
  environment}. We adopt here the standard ``quenched'' point of view:
we assume a realization of the sequence $(G_n)_{n \geq 2}$ is given,
and we aim at proving results that hold for almost every realization
of the graph sequence. \\ We consider the following random evolution,
indexed by the same time variable $n \geq 2$.  Let $S$ be a finite
set, representing possible choices made by ``individuals'' $i \in
V_n$.  To each vertex~$i \in V_n$ is associated a probability~$p_{n,i}
\in {\mathcal P}(S)$.  The quantity $p_{n,i}(x)$ ($x \in S$)
represents the inclination of individual~$i$ to adopt the choice~$x$
at time~$n$ or, in different terms, the {\em relative opinion} of
individual~$i$ about~$x$: the higher this value, the better the
opinion of~$i$ on $x$ compared with that on the other alternatives $y
\neq x$ ($y \in S$).  The following two-steps dynamics occurs before
the arrival of the $(n+1)^{\mbox{\tiny{st}}}$ vertex.
\medskip
\bi
\item[{\it Step~1}\,:] ({\em consensus}) Through a {\em fast}
  consensus dynamics on the graph $G_n$, the $p_{n,i}$ are
  homogenized: every vertex ends up with the same inclination
$$ \mathcal{Z}_{n,i} = \mathcal{Z}_{n}:=\frac{1}{n} \sum_{i=1}^n p_{n,i}.$$

\item[{\it Step~2}\,:] ({\em hub's influence}) Let $j_n$ be a vertex
  chosen arbitrarily among those of maximal degree.  This vertex
  exhibits a choice $I_{n+1}=x$ with probability $\mathcal{Z}_n(x)$,
  i.e. according to his (and everyone else's) inclination. The
  exhibition of the choice has influence on the inclination of the
  vertex~$j_n$'s neighbors, so that, given $I_{n+1}=x$:
\[
\begin{array}{ll}
 p_{n+1,j}=\l \delta_x +(1-\l) \mathcal{Z}_n &
\mbox{if $j$ is a neighbor of $j_n$} \\
 p_{n+1,j}= \mathcal{Z}_n & \mbox{otherwise,}
 \end{array}
 \]
 where $\l \in (0,1)$ is a given constant and $\d_x$ denotes the Dirac
 measure at $x \in S$.  \ei After these two steps, the vertex $n+1$ is
 added; its inclination $p_{n+1,n+1}$ right after arrival could be
 taken arbitrarily; just for simplicity in next formulas, we set
 $p_{n+1,n+1} = \mathcal{Z}_n$.

 \noindent This dynamics allows to obtain a recursive formula for
 $\mathcal{Z}_n$:
 \[
 \begin{split}
\mathcal{Z}_{n+1} & =\frac{1}{n+1}
\left[ \left(n +1-d_{j_n}(n)\right) \mathcal{Z}_n
+ d_{j_n}(n) (\l \delta_{I_{n+1}} + (1-\l) \mathcal{Z}_n) \right] \\
& = \frac{1}{n+1} \left[ (n+1 -\l \, d_{j_n}(n))\mathcal{Z}_n +
\l d_{j_n}(n) \delta_{I_{n+1}} \right]  \\
&  = (1-r_n) \mathcal{Z}_n + r_n \delta_{I_{n+1}},
\end{split}
\]
where
\be{rn-graph}
r_n := \frac{\l d_{j_n}(n)}{n+1}.
\ee
This has the form
\eqref{dyn-symplex} with $K_n(y) = \d_y$.\\
\smallskip

It should be stressed that many variants of this consensus-influence
dynamics could be considered as well; for instance, influence could be
exercised by vertices other than those with maximal degree, e.g. with
a degree dependent probability. Our specific choice makes particularly
easy to verify the conditions of some of the results below, see Remark
\ref{rem:examples} for details.
\\
\bigskip

\indent Finally, we remark that the dynamics \eqref{dyn-symplex} are
special cases of \emph{stochastic algorithms}, that are treated with
stochastic approximation methods and are used in many different
contexts and applications (see \cite{Bena} for an overview,
\cite{Kush} for a general reference and applications and \cite{fabian}
for classical results in the spirit of this paper). In particular,
\eqref{dyn-symplex} admits the following algorithmic
interpretation. Let $K: S \ra \mathcal{P}(S)$ be given. It can be
viewed as a stochastic kernel that induces a map $T_K: \mathcal{P}(S)
\ra \mathcal{P}(S)$ by
\[
T_K \mu := \sum_{y \in S} \mu(y)K(y).
\]
Then, \eqref{dyn-symplex} is a version of the Robbins-Monro algorithm
(see \cite{Kush}) to obtain a fixed point of $T_K$, i.e. a stationary
distribution of the $S$-valued Markov chain with transition kernel
$K$.  \\

\bigskip
\indent This paper is concerned with systems of $N$ interacting random
walks in which, to $N$ evolutions as in \eqref{dyn-symplex}, we add an
{\em interaction} term of {\em mean-field} type.  We are particularly
interested in the phenomenon of {\em synchronization}, that could be
roughly defined as the tendency of different components to adopt a
common long-time behavior. This phenomenon has been subject to recent
investigation in systems of many interacting particles, where
synchronization emerges in the large-scale limit \cite{daipra15, GB,
  perth1}. More recently, interacting urn models have attracted
attention as prototypical dynamics subject to reinforcement
\cite{aletti, ben, marsili, paganoni, pemantle-2007}. For some of
these dynamics, synchronization is induced by reinforcement, so it
does not require a large-scale limit \cite{cri-dai-min, daipra-2014,
  launay-2012, launay}. Another context in which synchronization
emerges naturally is that of opinion dynamics in a population
\cite{fagnani}. We have proposed here in Example~1.2 a version of
opinion dynamics in an evolving population: the interacting version
could be interpreted as related to different homogeneous groups within
a given population, in the same spirit as in \cite{collet, contucci}.
\\
\bigskip

\indent In this work, to avoid complications, we focus on the case $S
= \{0,1\}$, so that there is only one relevant variable, $Z_n :=
\mathcal{Z}_n(1)$; moreover we assume $K_n(y) = K(y)$ to be
independent of time. Concerning the examples considered above, this
time-independence property holds for the opinion models in
preferential attachment graphs; in urn models a sufficient condition
is that the reinforcement matrix $A_n$ is of the form $A_n = c_n A$
for some $c_n >0$ and a given matrix $A$ independent of $n$. This
includes generalizations of P\'olya and Friedman models, where the
reinforcement matrix is allowed to depend on $n$. The most general
function $K: S \ra \mathcal{P}(S)$ can be written in the form
\be{model0} K(y) = \rho \d_y + (1-\rho) q, \ee for some $\rho \in
   [0,1]$ and a given $q \in \mathcal{P}(S)$. After identifying $q$
   with $q(1)$, the evolution of the $i$-th walk is therefore given by
   \be{model1} Z_{n+1}(i) = (1-r_n) Z_n(i) + r_n (\rho I_{n+1}(i) +
   (1-\rho)q).  \ee

The interaction enters in the conditional law of $I_{n+1}(i)$: Setting
\[
\F_n := \s(Z_k(i): \, i = 1,2,\ldots,N\,; \; 0 \leq k \leq n)\,,
\]

we assume that the random variables $\{I_{n+1}(i), \, i =
1,2,\ldots,N\}$ are conditionally independent given $\F_n$ with

\be{model2} P(I_{n+1}(i) = 1 | \F_n) = (1-\alpha)  Z_n(i) + \alpha  Z_n
\ee

where $\alpha \in [0,1]$ is the interaction parameter and
\be{model3}
Z_n := \frac{1}{N} \sum_{i=1}^N Z_n(i).
\ee

Under suitable conditions on $r_n$, but actually no conditions if
$\rho = 1$, the sequence $(Z_n)$ converges almost surely to a limit
$Z$. In the spirit of similar results for urn models, we study the
corresponding rate of convergence and, in some cases, we obtain a
fluctuation Theorem in functional form. We compare this rate of
convergence with the {\em rate of synchronization}, which we define as
the rate at which $Z_n(i) - Z_n$ converges to zero. As observed in
\cite{neeraja} for interacting Friedman urns, synchronization may be
     {\em faster} than convergence. In our model we show that this
     occur when $\rho=1$ and $r_n \sim \frac{c}{n^{\g}}$ with $c>0$
     and $\frac12 < \g < 1$. We stress the fact that this is not a
     large-scale phenomenon, in the sense that it holds {\em for any}
     value of $N$.  \\
\bigskip

\indent The paper is organized as follows. In section 2 we present our
main results. Section 3 contains some basic identities often used in
the proofs. Sections 4-7 are then devoted to proofs.

\section{Main results}

From now on we study interacting dynamics of the form defined by
\eqref{model1}, \eqref{model2} and \eqref{model3}. We assume that the
initial configuration $[Z_0(i)]_{i=1}^N$ has a permutation invariant
distribution with $\E[Z_0(i)] = \frac12$ and $\E[Z_0(1-Z_0)]>0$ where
$Z_0:=N^{-1}\sum_{i=1}^N Z_0(i)$. The value $\frac12$ could be
replaced by $z_0 \in (0,1)$ at the only cost of longer formulas. It is
worthwhile also to note that these assumptions will be used only in some
of our proofs and they could be weakened.

\subsection{Convergence and synchronization}

The following theorem describes the convergence of the sequence $(Z_n)$.

\bt{prop:conv1}
\mbox{}\\ \vspace{-4mm}
\bi
\item[(i)] If $\rho =1$, then $Z_n$ converges almost surely to a random
  variable $Z$. Moreover:
\begin{itemize}
\item[a)] If $\alpha>0$, then
\be{conv1}
\P(Z \in \{0,1\}) = 1 \Leftrightarrow
\sum_{n} r_n^2 = +\infty;
\ee
\item[b)] If $\alpha=0$, each $Z_n(i)$ converges almost surely to a
  random variable $Z(i)$ such that
\be{conv1-irene}
\P(Z(i) \in \{0,1\}) = 1 \Leftrightarrow \sum_{n} r_n^2 =
+\infty.
\ee
\end{itemize}

\item[(ii)]
If $\rho <1$ and
\be{conv2bis}
\sum_{n} r_n= +\infty \ \ \ \mbox{and} \ \ \
\sum_{n} r_n^2 < +\infty\,,
\ee
then $Z_n \ra q$ almost surely.
\ei
\et

The following result particularly points out that in the case of a
single walk, strengthening condition \eqref{conv1} one gets the
phenomenon of {\em fixation}. Note that this phenomenon has been
observed in various urn models, see e.g. \cite{davis,launay-2012,
  launay}.

\bp{prop:fixation}
Assume $N=1$. If\\

\be{conv2}
\rho=1\quad\hbox{ and}\quad  \sum_{n \geq 1} \exp \left[ -
  \sum_{k=0}^n r_k^2 \right] < +\infty
\ee
or\\

\be{conv2-bis}
\rho<1,\quad q\in\{0,1\} \quad\hbox{ and}\quad
\sum_n r_n=+\infty\,,
\ee
then there exists a random index $M$ such that with probability one
the indicator functions $\{I_n: \, n \geq M\}$ have all the same
value.
\ep

Note that, if $r_n = O(n^{-\g})$, then \eqref{conv1} holds for $\g
\leq \frac12$, while \eqref{conv2} for $\g<\frac12$. Moreover,
\eqref{conv2bis} holds for $\frac12 < \g \leq 1$.  \\

Next theorem establishes the fact that {\em synchronization} indeed takes
place as soon as either interaction is present ($\alpha>0$) or the limit
of $Z_n$ is deterministic ($\rho<1$).

\bt{prop:conv2}
Suppose that \eqref{conv2bis} holds and $\rho(1-\alpha) <1$. 
Then, for all $i \in \{1,2,\ldots,N\}$, we have 
\[
Z_n(i) - Z_n \longrightarrow 0 \ \ \ \mbox{a.s.}
\]
In particular, if $Z$ is the almost sure limit of $Z_n = \frac{1}{N}
\sum_{i=1}^N Z_n(i)$ (note that, for $\rho<1$, by Theorem
\ref{prop:conv1}, $Z=q$), we have $Z_n(i) \ra Z$ almost surely.  
\et

As we will see, the proof of this result does not require the
assumptions on the initial configuration.

\br{rem:examples} At this point it is worth discussing the assumptions
in the previous results, in the context of the applications 1.1 and
1.2 proposed in the introduction. \\

{\em Urn models}. Consider an urn model with reinforcement matrix
$A_n$ of the form $A_n = c_n A$. So, after the $(n+1)$-st drawing, the
number of balls added into the urn is $\overline{A}_n = c_n
\overline{A}$ and therefore
\[
r_n = \frac{c_n \overline{A}}{N(0) + \overline{A} \sum_{k=0}^n c_k}.
\]

As observed in the introduction, $r_n$ is of order $\frac{1}{n}$
whenever $c_n$ grows polynomially, so \eqref{conv2bis} holds. For a
different behavior one has to consider a faster growing reinforcement,
e.g $c_n = \exp\left(n^{\b}\right)$. In this case, for $0<\b<1$, it is
easily shown that $r_n$ is of order $\frac{1}{n^{1-\b}}$, thus
\eqref{conv1}, \eqref{conv2} or \eqref{conv2bis} may hold depending on
the value of $\b$. \\

{\em Opinion dynamics}. Note that in this model $\rho = 1$. By Theorem
8.8 in \cite{remco}, the maximal degree $d_{j_n}(n)$ at time $n$ is
such that the limit
\[
\lim_{n \ra +\infty} \frac{d_{j_n}(n)}{n^{\frac{1}{2+\d}}} =: l
\]
exists for almost every realization of the graph sequence $(G_n)_{n
  \geq 2}$. Thus, given the definition \eqref{rn-graph} of $r_n$, we
have $r_n \sim \frac{c}{n^{\g}}$, with $\g = \frac{1+\d}{2+\d}$ and $c
= \l l$.  It follows that the conditions \eqref{conv2bis} hold for
every $\d>0$. For $-1 < \d \leq 0$, condition \eqref{conv1} holds and
so the population's inclination ``polarizes'', i.e. it converges to
the Dirac measure concentrated on one choice.  \er

\subsection{Fluctuation theorems}

Assume $\rho(1-\alpha) <1$ and $r_n \sim \frac{c}{n^{\g}}$, where $\frac{1}{2}
< \g \leq 1$, with the meaning
\[
\lim_{n \ra +\infty} n^{\g} r_n = c > 0.
\]

Note that the above assumptions imply that Theorem \ref{prop:conv2}
holds. In all the following theorems, the notation
$\stackrel{d}\longrightarrow$ denotes convergence in distribution with
respect to the classical Skorohod's topology (see e.g. \cite{Bill2}).

\bigskip

Next result describes the fluctuations of $Z_n$ around its limit
$Z$ in terms of a functional Central Limit Theorem in the case $\rho
= 1$.

\bt{th:fluctuationsZ}
Suppose $\rho = 1$ (and so $\alpha>0$). Then the random limit $Z$ of
$Z_n$ is such that $\P(Z \in \{0,1\}) < 1$ and $\P(Z = z) = 0$ for all
$z \in (0,1)$.  Moreover the following holds:

\begin{equation}\label{fcltZ}
\left( t^{2\gamma -1} n^{\g- \frac12} (Z_{\lfloor n t\rfloor}-Z)
\right)_{t\geq 0}
\stackrel{d}\longrightarrow (W_{V_t})_{t\geq 0}
\end{equation}

where
$$
V_t=\frac{c^2}{N(2\gamma-1)} Z(1-Z) t^{2\gamma-1}
$$
and $W=(W_t)_{t\geq 0}$ is a Wiener process independent of $V=(V_t)_{t\geq 0}$.
\et

Next theorem characterizes the rate of synchronization, i.e. the rate
of convergence to zero of $Z_n(i) - Z_n$, in terms of a functional
Central Limit Theorem in the case $\rho = 1$.
\\

\bt{FCLTSyncro}
Suppose $\rho = 1$ (and so $\alpha>0$).
\bi
\item[(i)]
If $\frac12 < \g < 1$ then
\[
\left(n^{\g/2} e^{c\alpha t}
\left(Z_{\lfloor n+n^{\g}t \rfloor}(i) - Z_{\lfloor n+n^{\g}t \rfloor} \right)
\right)_{t \geq 0}
\stackrel{d}\longrightarrow (W_{V_t})_{t \geq 0},
\]
where
\[
V_t = \left(1-\frac{1}{N}\right)\frac{c Z(1-Z)}{2\alpha } e^{2c\alpha  t}.
\]
and $W=(W_t)_{t \geq 0}$ is a standard Brownian motion independent of
$V=(V_t)_{t\geq 0}$.

\item[(ii)]
If $\g = 1$ and $2c\alpha  >1$, then
\[
\left(n^{1/2} (1+t)^{c\alpha }
\left(Z_{\lfloor n + nt \rfloor}(i) - Z_{\lfloor n+ nt \rfloor} \right)
\right)_{t \geq 0}
\stackrel{d}\longrightarrow (W_{V_t})_{t \geq 0},
\]
where
\[
V_t = \left(1-\frac{1}{N}\right)
\frac{ c^2Z(1-Z)}{2c\alpha -1} (1+t)^{2c\alpha -1}.
\]
\ei
and $W=(W_t)_{t \geq 0}$ is a standard Brownian motion independent
of $V=(V_t)_{t\geq 0}$
\et

\br{rem:synchro} When $\rho=1$, since $\gamma/2 > \gamma -1/2$ for
$\frac12 < \g < 1$, by Theorems \ref{th:fluctuationsZ} and
\ref{FCLTSyncro}, we have that in this regime synchronization is
faster than convergence. More precisely, the proof of Theorem
\ref{FCLTSyncro} implicitly contains the fact that, for
$1/2<\gamma<1$,
\begin{equation}\label{nuova-*}
\E\left[(Z_n(i)-Z_n)^2\right]\sim \left(1-\frac{1}{N}\right)C_1\, n^{-\gamma}
\end{equation}
with a suitable constant $C_1 = c (2\alpha)^{-1} \E[Z(1-Z)]>0$ by
Theorem \ref{prop:conv1}. This fact, by permutation invariance,
implies that we have for $i\neq j$
\begin{equation*}
\E\left[(Z_n(i)-Z_n(j))^2\right]=2 \frac{N}{N-1}\E\left[(Z_n(i)-Z_n)^2\right]
\sim 2C_1\, n^{-\gamma}.
\end{equation*}
On the other hand, the proof of Theorem \ref{th:fluctuationsZ} (see
also the proof of Proposition \ref{th:stable}) implicitly contains the
fact that, for $1/2<\gamma<1$,
\begin{equation*}
\E\left[(Z_n-Z)^2\right]\sim \frac{1}{N}C_2\, n^{-(2\gamma-1)}
\end{equation*}
with a suitable constant $C_2 = c^2 (2\gamma-1)^{-1} \E[Z(1-Z)]>0$ by
Theorem \ref{prop:conv1}. This fact, together with \eqref{nuova-*},
implies that, for each $i$, we have
\begin{equation*}
\E\left[(Z_n(i)-Z)^2\right]\sim \E\left[(Z_n-Z)^2\right]
\sim \frac{1}{N}C_2\, n^{-(2\gamma-1)}.
\end{equation*}
Therefore, for $\rho=1$ and $1/2<\gamma<1$, the velocity of
convergence to zero of $\|Z_n(i)-Z_n(j)\|_{L^2}$ is greater than the
one of $\|Z_n(i)-Z\|_{L^2}$.
\er

In the case of $\rho <1$ and $q \not\in \{0,1\}$ the rate of
convergence of $Z_n$ to its limit $q$ is, for $\frac12<\g<1$,
different from the scaling in Theorem \ref{th:fluctuationsZ}, and
matches that in Theorem \ref{FCLTSyncro}, that is $\g/2$. The
following two results, in particular, show that convergence and
synchronization occur at the same rate.

\bt{FCLT-q}
Suppose $\rho < 1$ and $q \not\in \{0,1\}$.
\begin{itemize}
\item[(i)]
If $\frac12 < \g < 1$, then
\[
\left(n^{\g/2} e^{c(1-\rho)t}
\left(Z_{\lfloor n+n^{\g}t \rfloor} - q \right)
\right)_{t \geq 0}
\stackrel{d}\longrightarrow (W_{V_t})_{t \geq 0},
\]
where $W=(W_t)_{t \geq 0}$ is a standard Brownian motion and
\[
V_t =\frac{c q(1-q)\rho^2}{2N(1-\rho)} e^{2c(1-\rho) t}.
\]

\item[(ii)]
If $\g = 1$ and $2c(1-\rho) >1$, then
\[
\left(n^{1/2} (1+t)^{c(1-\rho)}
\left(Z_{\lfloor n+nt \rfloor} - q \right)
\right)_{t \geq 0}
\stackrel{d}\longrightarrow (W_{V_t})_{t \geq 0},
\]
where $W=(W_t)_{t \geq 0}$ is a standard Brownian motion and
\[
V_t = \frac{ c^2 q(1-q)}{N(2c(1-\rho)-1)}(1+t)^{2c(1-\rho)-1}.
\]
\end{itemize}
\et

\bt{FCLTSyncro-q}
Suppose $\rho < 1$ and $q \not\in \{0,1\}$.
\begin{itemize}
\item[(i)]
If $\frac12 < \g < 1$, then
\[
\left(n^{\g/2} e^{c(1-\rho(1-\alpha) )t}
\left(Z_{\lfloor n+n^{\g}t \rfloor}(i) - Z_{\lfloor n+n^{\g}t \rfloor} \right)
\right)_{t \geq 0}
\stackrel{d}\longrightarrow (W_{V_t})_{t \geq 0},
\]
where $W=(W_t)_{t \geq 0}$ is a standard Brownian motion and
\[
V_t =
\left(1-\frac{1}{N}\right)\frac{c\rho^2 q(1-q)}{2(1-\rho(1-\alpha) )}
e^{2c(1-\rho(1-\alpha) ) t}.
\]

\item[(ii)]
If $\g = 1$ and $2c(1-\rho(1-\alpha) ) >1$, then
\[
\left(n^{1/2} (1+t)^{c(1-\rho(1-\alpha) )}
\left(Z_{\lfloor n+ nt \rfloor}(i) - Z_{\lfloor n+ nt \rfloor}\right)
\right)_{t \geq 0}
\stackrel{d}\longrightarrow (W_{V_t})_{t \geq 0},
\]
where $W=(W_t)_{t \geq 0}$ is a standard Brownian motion and
\[
V_t =
\left(1-\frac{1}{N}\right)
\frac{ c^2 q(1-q)}{2c(1-\rho(1-\alpha) )-1}(1+t)^{2c(1-\rho(1-\alpha) )-1}.
\]
\end{itemize}

\et

\br{comparison} Functional Central Limit Theorems in the spirit of
those above have been proved for various urn models
(e.g. \cite{bai-et-al,bas-das, gouet, janson, z}). In particular
\cite{bai-et-al, gouet} and \cite{z} contain results for Friedman urn
models and P\'olya urn models respectively, that in our model
correspond to the case $N=1$ and $\g=1$.  The results for the
fluctuations of the Friedman urn, in particular, show that the
condition $2c(1-\rho(1-\alpha) ) >1$ in Theorem \ref{FCLTSyncro-q} is
essential: the Friedman urn that corresponds to $2c(1-\rho(1-\alpha) )
<1$ is known to have non-Gaussian fluctuations (see \cite{fr65}), so
no convergence to a Gaussian process is possible. The case of one
Friedman urn with $2c(1-\rho) = 1$ is considered in \cite{gouet}: a
functional central limit theorem holds with a logarithmic correction
in the scaling. We do not consider this case here.\\ For interacting
P\'olya urns ($\g=1$) a non functional version of Theorem
\ref{FCLTSyncro} is proved in \cite{cri-dai-min}, under the same
condition $2c\alpha >1$. For interacting Friedman urns ($\g=1$) a non
functional version of Theorem \ref{FCLTSyncro-q} is proved in
\cite{neeraja}.\\ We finally remark that the non-functional version of
Theorem \ref{FCLT-q} could be alternatively derived by following
stochastic approximation methods (see e.g. \cite{Bena, fabian, Kush}).
\er

\br{casispeciali} Regarding the assumption on $q$ of Theorems
\ref{FCLT-q} and \ref{FCLTSyncro-q}, we note that, if $q\in\{0,1\}$,
when $\gamma<1$, the behaviors of $Z_n$ and $Z_n(i)$ are ``eventually
deterministic'' (as can be easily derived by the same argument used
for Proposition \ref{prop:fixation}). Therefore the only case to be
considered is when $\g = 1$, but we will not deal with it in this
paper. For $N=1$, a functional central limit Theorem could be
obtained as in Proposition 2.2 of \cite{gouet}.
\er

\section{Basic properties}

In this section we derive some simple recursions on the random walks
$Z_n(i)$ and on $Z_n=\frac{1}{N}\sum_{i=1}^N Z_n(i)$, that will be
used several times.
\\

By averaging over $i$ in \eqref{model1}, we have
\[
Z_{n+1}-Z_n=r_n\left[\rho \left(N^{-1}\sum_{i=1}^N I_{n+1}(i)-Z_n\right) -
(1-\rho) (Z_n-q) \right]
\]
with
\[
\mathrm{E}\left[N^{-1}\sum_{i=1}^N I_{n+1}(i)\Big|\F_n\right]=
Z_n.
\]

Therefore we can write
\begin{equation}\label{eq-base-inter}
Z_{n+1}=Z_n - r_n (1-\rho) (Z_n-q) + \rho r_n \Delta M_{n+1},
\end{equation}
where
\begin{equation}\label{def-M_n}
\Delta M_{n+1}=
N^{-1}\sum_{i=1}^N I_{n+1}(i)-
\mathrm{E}\left[N^{-1}\sum_{i=1}^N I_{n+1}(i)\Big|{\cal F}_n\right]
=N^{-1}\sum_{i=1}^N I_{n+1}(i)-Z_n.
\end{equation}

Then, subtracting  \eqref{eq-base-inter} to   \eqref{model1}, we obtain
\be{eq-base-diff}
Z_{n+1}(i) - Z_{n+1} =
[1-r_n (1- \rho (1-\alpha) )] (Z_n(i) - Z_n) +
r_n \rho [\Delta M_{n+1}(i) - \Delta M_{n+1}],
\ee
where
\be{def-Mi}
\Delta M_{n+1}(i) =I_{n+1}(i) - \E[I_{n+1}(i) | \F_n].
\ee

In particular, from the relations above, we have
\begin{eqnarray}
\mathrm{E}[Z_{n+1}-q|\mathcal{F}_n]&=& [1-(1-\rho)r_n](Z_n-q)\label{s_cond_zn-q}
\\
\mathrm{E}[Z_{n+1}(i)-Z_{n+1}|\mathcal{F}_n]&=&
[1-(1-\rho(1-\alpha) )r_n]\left(Z_n(i)-Z_n\right)\label{s_cond_zni-zn}
\end{eqnarray}
and
\begin{eqnarray}
\mathrm{Var}[Z_{n+1}(i)-Z_{n+1}|\mathcal{F}_n]&=&
r_n^2\rho^2
\mathrm{E}\left[( \Delta M_{n+1}(i)-\Delta M_{n+1} )^2 | {\cal F}_n\right]
\nonumber
\\
&=&
r_n^2\rho^2\left(1-\frac{1}{N}\right)^2\mathrm{Var}[I_{n+1}(i)|\mathcal{F}_n]
+\frac{r_n^2\rho^2}{N^2}\sum_{j\neq i}\mathrm{Var}[I_{n+1}(j)|\mathcal{F}_n]
\label{ric-var1}
\\
\mathrm{Var}[Z_{n+1}|\mathcal{F}_n]&=&\frac{r_n^2\rho^2}{N^2}
\sum_{j=1}^N\mathrm{Var}[I_{n+1}(j)|\mathcal{F}_n].   \label{ric-var2}
\end{eqnarray}

\section{Proofs: convergence and synchronization}

\subsection{Proof of Theorem \ref{prop:conv1}}
\mbox{}\\
\noindent{ \em Part (i)(a)} Here we assume $\rho =1$ and
$\alpha>0$. By \eqref{eq-base-inter}, we immediately get that $(Z_n)$
is a bounded martingale. Therefore, it converges a.s. (and in $L^p$)
to a random variable $Z$, with values in $[0,1]$. \\

\indent Since by assumption $\E[Z_0(i)] = \frac12$ for every $i$, we
have $\E(Z) = \frac12$.  Moreover $P(Z \in \{0,1\}) = 1$ if and only
if
\[
\var(Z) = \lim_{n \ra +\infty} \var(Z_n) = \frac14.
\]
By using \eqref{ric-var2}, we have 
\begin{eqnarray*}
\mathrm{Var}[Z_{n+1}]&=&\mathrm{E}[\mathrm{Var}(Z_{n+1}|\mathcal{F}_n)]+
\mathrm{Var}[\mathrm{E}(Z_{n+1}|\mathcal{F}_n)]\\
&=& r_n^2
N^{-2} \sum_{i=1}^N
\mathrm{E}\left[
\left( (1-\alpha)  Z_n(i)+\alpha  Z_n \right)
\left(1-(1-\alpha)  Z_n(i)-\alpha  Z_n \right)
\right]
+\mathrm{Var}[Z_n]
\\
&=&
r_n^2
\mathrm{E}\left[
Z_n/N-(\alpha^2+2(1-\alpha) \alpha )Z_n^2/N
-(1-\alpha)^2\sum_{i=1}^N Z_n^2(i)/N^2
\right]
+\mathrm{Var}[Z_n]
\\
&=&
r_n^2
\mathrm{E}\left[
Z_n/N-(1-(1-\alpha)^2)Z_n^2/N-(1-\alpha) ^2Z_n^2
+(1-\alpha) ^2\sum_{i\neq j}Z_n(i)Z_n(j)/N^2
\right]
+\mathrm{Var}[Z_n]
\\
&=&
[1/N - (1-\alpha)^2 (1-1/N)] r_n^2/4
+(1-\alpha) ^2 r_n^2 \sum_{i\neq j} \mathrm{E}[Z_n(i)Z_n(j)]/N^2
\\
&+&\left\{1 - [1/N + (1-\alpha) ^2(1-1/N)] r_n^2 \right\} \mathrm{Var}[Z_n]
\\
&=&
[1/N + (1-\alpha)^2 (1-1/N)] r_n^2/4 - (1-\alpha)^2(1-1/N)r_n^2/2
+(1-\alpha)^2 r_n^2 \sum_{i\neq j}\mathrm{E}[Z_n(i)Z_n(j)]/N^2
\\
&+&\left\{1 - [1/N + (1-\alpha)^2(1-1/N)] r_n^2 \right\}\mathrm{Var}[Z_n].
\end{eqnarray*}
Now we observe that, using the permutation invariance, we have
$\E[Z^2_n(i)] \leq \E[Z_n(i)] = \frac12$ for all $i$ and
\begin{equation*}
\sum_{i\neq j}\mathrm{E}[Z_n(i)Z_n(j)]/N^2 =  \E[Z^2_n] -
\frac{1}{N}\E[Z^2_n(1)] 
 \geq  \var[Z_n] + \frac14 - \frac{1}{2N}.
\end{equation*}
Using this fact in the formula above, we get
\[
\var[Z_{n+1}]  \geq \left[1- \frac{1}{N}(1-(1-\alpha) ^2) r^2_n \right]
\var[Z_n] +
\frac{1}{N}(1-(1-\alpha) ^2) \frac{r^2_n}{4},
\]
which, by letting $x_n := \frac14 - \var[Z_n] \geq 0$, is equivalent to
\be{polar}
x_{n+1} \leq  \left(1-C r^2_n \right) x_n
\ee
with $C:= \frac{1-(1-\alpha) ^2}{N}\in\, ]0,1]$ for $0\leq (1-\alpha)
    <1$. Therefore
\[
x_n \leq x_0 \prod_{k=0}^{n-1} \left(1-C r^2_k \right)
\]
which implies $x_n \ra 0$ if $\sum_n r_n^2 = +\infty$.\\

\indent We are left to prove that if $\sum_n r_n^2 < +\infty$ then
$x_n \not\rightarrow 0$. From the above equalities, we have
\begin{eqnarray*}
\mathrm{Var}[Z_{n+1}]&=&
 r_n^2
N^{-2} \sum_{i=1}^N
\mathrm{E}\left[
(1-\alpha)  Z_n(i)+\alpha  Z_n 
\right]
-
\mathrm{E}\left[
\left( (1-\alpha)Z_n(i)+\alpha  Z_n \right)^2
\right]
+\mathrm{Var}[Z_n]
\\
&=&
\frac{r_n^2}{N} \mathrm{E}[Z_n]
-
\frac{r_n^2}{N^2}
\sum_{i=1}^N \mathrm{E}\left[
\left( Z_n+(1-\alpha)(Z_n(i)-Z_n) \right)^2
\right]
+\mathrm{Var}[Z_n]
\\
&=&
\frac{r_n^2}{N} \mathrm{E}[Z_n]
-
\frac{r_n^2}{N}\mathrm{E}[Z_n^2]
- 
\frac{r_n^2}{N^2}(1-\alpha)^2\sum_{i=1}^N\mathrm{E}\left[(Z_n(i)-Z_n)^2\right]
+\mathrm{Var}[Z_n]
\\
&\leq&
\frac{r_n^2}{N} \mathrm{E}[Z_n]
-
\frac{r_n^2}{N}\mathrm{E}[Z_n^2]
+\mathrm{Var}[Z_n]
=
\frac{r_n^2}{4N}+ \left(1-\frac{r_n^2}{N}\right)\mathrm{Var}[Z_n],
\end{eqnarray*}
where we have used the identities $\E(Z_n) = \frac{1}{2}$ and
$\E[Z_n^2] = \frac{1}{4} + \var[Z_n]$. Thus we have
\begin{equation*}
x_{n+1}\geq \left(1-\frac{r_n^2}{N}\right)x_n, 
\end{equation*}
from which  it follows
\[
x_n \geq x_0 \prod_{k=0}^{n-1} \left(1-\frac{r_k^2}{N}\right),
\]
where $x_0>0$ by assumption. Since, $\sum_{n} r_n^2 < +\infty$ by
assumption, we obtain $\lim_{n \ra +\infty} x_n > 0$.  \\
\smallskip

\noindent{\em Part(i)(b)}. In the case when $\rho=1$ and $\alpha=0$,
each $(Z_n(i))$ is a bounded martingale and so we have the almost sure
(and in $L^p$) convergence of $Z_n(i)$ to a random variable $Z(i)$ with
values in $[0,1]$. Moreover, with similar computation as above, we have
$$
\var[Z_{n+1}(i)]=\frac{r_n^2}{4}+(1-r_n^2)\var[Z_n(i)]
$$ that is \be{ric-new} x_{n+1}(i)=(1-r_n^2)x_n(i) \ee where
$x_n(i):=\frac{1}{4}-\var[Z_n(i)]$. Therefore the conclusion
immediately follows.\\
\smallskip

\noindent
{\em Part (ii)}. We are now assuming $\rho <1$ and \eqref{conv2bis},
that is $ \sum_{n} r_n = +\infty$ and $\sum_{n} r_n^2 < +\infty$.
\\

Using \eqref{eq-base-inter}, we have \be{proof3} {\mathrm
  E}\left[\left(Z_{n+1}-q\right)^2|{\cal F}_n\right] =
(Z_n-q)^2[1-2(1-\rho) r_n] +r_n^2 \left\{ (1-\rho)^2 (Z_n-q)^2 +
\rho^2{\mathrm E}\left[(\Delta M_{n+1})^2|{\cal F}_n\right] \right\}
\ee from which we obtain
\[
{\mathrm E}\left[\left(Z_{n+1}-q\right)^2|{\cal F}_n\right]
\leq (Z_n-q)^2+ r_n^2 \xi_n
\]
where $\xi_n:=\left\{ (1-\rho)^2 (Z_n-q)^2 + \rho^2{\mathrm
  E}\left[(\Delta M_{n+1})^2|{\cal F}_n\right]
\right\}$ is bounded.
Therefore, since
$\sum_n r_n^2 <+\infty$, we can conclude that
$\big((Z_n-q)^2\big)_n$ is a positive almost supermartingale (see
\cite{rob-sie}) and so it converges almost surely (and in $L^p$).
In order to show that this limit is $0$, we are left to show that
\be{proof4}
\lim_{n \ra +\infty}
\E\left[(Z_n - q)^2 \right] = 0.
\ee
Averaging in \eqref{proof3} and letting
$x_n := {\mathrm E}\left[\left(Z_{n}-q\right)^2\right]$, we have
\[
x_{n+1} = [1-2(1-\rho) r_n] x_n + K_n r_n^2,
\]
with
\[
0 \leq K_n := \left\{
(1-\rho)^2 \E[(Z_n-q)^2] +
\rho^2{\mathrm E}\left[(\Delta M_{n+1})^2\right]
\right\} \leq 1.
\]
Recalling the assumptions \eqref{conv2bis} and $\rho<1$, the
conclusion follows from Lemma \ref{lemma-tecnico}.
\qed

\subsection{Proof of Proposition \ref{prop:fixation}}

First suppose that \eqref{conv2} holds.  Setting
$x_n:=\frac{1}{4}-\var[Z_n]$ and using \eqref{ric-new}, we have, for a
suitable positive constant $C$,

\be{stimavar} x_{n} = x_0 \prod_{k=0}^{n-1}
(1-r_k^2) \sim C \exp \left[ - \sum_{k=0}^n r_k^2 \right].  \ee Now
\begin{eqnarray*}
\sum_{n=1}^{+\infty}P(I_{n}=0, I_{n+1}=1)&=&
\sum_{n=1}^{+\infty}\mathrm{E}[P(I_{n}=0, I_{n+1}=1|\mathcal{F}_{n})]
= \sum_{n=1}^{+\infty}\mathrm{E}[(1-I_n)Z_n]\\
&=& \sum_{n=1}^{+\infty}\mathrm{E}[(1-I_n)(Z_{n-1}(1-r_{n-1})+r_{n-1}I_n)]\\
&=& \sum_{n=1}^{+\infty}\mathrm{E}[(1-I_n)Z_{n-1}(1-r_{n-1})]
\;\mbox{ then, conditioning on } \mathcal{F}_{n-1},\\
&=&\sum_{n=1}^{+\infty}(1-r_{n-1})\mathrm{E}[Z_{n-1}(1-Z_{n-1})]\\
&\leq& \sum_{n=1}^{+\infty}\mathrm{E}[Z_{n-1}(1-Z_{n-1})]
=  \sum_{n=1}^{+\infty} x_{n-1}<+\infty
\end{eqnarray*}
by \eqref{stimavar} and \eqref{conv2}. Then by
Borel-Cantelli lemma $P(\limsup_{n}\{I_{n}=0,
I_{n+1}=1\})=0$ and the conclusion follows.
\\

Now, assume \eqref{conv2-bis} and $q=0$ (the case $q=1$ is specular),
by a similar argument as above, we get
\[
\sum_n P(I_{n+1}=1) = \sum_n \E(Z_n) < +\infty
\]
since $E[Z_n]\sim C\exp[-(1-\rho)\sum_{k=0}^n r_k]$.
\qed

\subsection{Proof of Theorem \ref{prop:conv2}}

We aim at showing that \be{for:conv} Z_n(i) - Z_n \longrightarrow 0
\ \ \mbox{a.s } \ee Set $x_n := \E\left[(Z_n(i) - Z_n)^2 \right]$. The
proof is essentially the same as that of Theorem \ref{prop:conv1},
(ii): we first show that $\left((Z_n(i) - Z_n)^2\right)$ is a positive
almost supermartingale, which implies almost sure (and in $L^p$)
convergence to a limit, and then we show that \be{limx_n} \lim_{n \ra
  +\infty} x_n = 0, \ee so that the limit of $Z_n(i) - Z_n$ is
a.s. zero.

By \eqref{eq-base-diff}, we obtain
\begin{multline} \label{syncf1}
\E[(Z_{n+1}(i) - Z_n)^2| \F_n] = [1-2 r_n (1-\rho(1-\alpha) )] (Z_n(i) - Z_n)^2
\\ + r_n^2 \left\{ (1-\rho(1-\alpha) )^2 (Z_n(i) - Z_n)^2 + \rho^2
\E\left[(\Delta M_{n+1}(i) - \Delta M_{n+1})^2 |\F_n \right]\right\}
\end{multline}
and so
\[
\E[(Z_{n+1}(i) - Z_n)^2| \F_n] \leq (Z_n(i) - Z_n)^2 + r_n^2 \xi_n,
\]
where $\xi_n := (1-\rho(1-\alpha) )^2 (Z_n(i) - Z_n)^2 + \rho^2\E\left[(\Delta
  M_{n+1}(i) - \Delta M_{n+1})^2 |\F_n \right]$ is bounded. Since
\mbox{$\sum_n r_n^2 <+\infty$}, this implies that $((Z_{n+1}(i) -
Z_n)^2)$ is a positive almost supermartingale.\\

It remains to prove \eqref{limx_n}. Taking the expected value in
\eqref{syncf1}, we obtain
\[
x_{n+1} = [1-2 (1-\rho(1-\alpha) )r_n] x_n + K_n r_n^2
\]
for a bounded sequence $(K_n)_n$ of positive numbers. Since we assume
\eqref{conv2bis} and $\rho(1-\alpha) <1$, the conclusion follows by
applying Lemma \ref{lemma-tecnico}.

\section{Fluctuation theorems I: proof of Theorem \ref{th:fluctuationsZ}}

The synchronization result in Theorem \ref{prop:conv2}
gives, for $\rho (1-\alpha) <1$,
\be{forlucca}
\mathrm{Var}[I_{n+1}(i)|\mathcal{F}_n]=\left((1-\alpha)
Z_n(i)+\alpha Z_n\right) \left(1-(1-\alpha)
Z_n(i)-\alpha Z_n\right)\longrightarrow Z(1-Z)  \ \mbox{a.s.}
\ee
 for all $i$. This,
together with \eqref{ric-var1} and \eqref{ric-var2}, implies
the following useful
relations:

\begin{eqnarray}
\mathrm{Var}[Z_{n+1}(i)-Z_{n+1}|\mathcal{F}_n]&\sim&
r_n^2\rho^2\left(1-\frac{1}{N}\right)Z(1-Z)\label{sim-var}\\
\mathrm{Var}[Z_{n+1}|\mathcal{F}_n]&\sim&
\frac{r_n^2\rho^2}{N}Z(1-Z). \label{sim-var-zn}
\end{eqnarray}

Before proving fluctuation theorems in the functional form, we prove a
fluctuation theorem in non-functional form, but with a stronger form
of convergence, the almost sure conditional convergence (see the
appendix for details). This result, which has independent interest, is
useful here to prove that the limit $Z$ has no point mass in
$(0,1)$. We recall that, in this section, we assume $\rho = 1$ and
$(1-\alpha) <1$.

\bp{th:stable}
Under the assumptions of Theorem \ref{th:fluctuationsZ},
$$ n^{\gamma-\frac12}(Z_n-Z)\stackrel{stably}\longrightarrow
\mathcal{N}\left(0, \frac{c^2}{N(2\gamma-1)}Z(1-Z)\right).$$

Moreover, the above convergence is in the sense of the almost sure
conditional convergence w.r.t. ${\cal F}=({\cal F}_n)$.
\ep

\bpr We want to apply Theorem \ref{fam_tri_vet_as_inf}. Let us
consider, for each $n\geq 1$, the filtration $({\cal
  F}_{n,h})_{h\in\N}$ and the process $(M_{n,h})_{h\in\N}$ defined by

\begin{equation*}
{\cal F}_{n,0}={\cal F}_{n,1}={\cal F}_n,\qquad
M_{n,0}=M_{n,1}=0
\end{equation*}

and, for $h\geq 2$,

\begin{equation*}
{\cal F}_{n,h}={\cal F}_{n+h-1},\qquad
M_{n,h}=n^{\gamma-\frac{1}{2}}(Z_n-Z_{n+h-1}).
\end{equation*}

By \eqref{eq-base-inter} (with $\rho=1$), it is easy to verify that,
with respect to $({\cal F}_{n,h})_{h\geq 0}$, the process
$(M_{n,h})_{h\geq 0}$ is a martingale which converges in $L^1$ (for
$h\to +\infty$) to the random
variable~$M_{n,\infty}:=n^{\gamma-\frac{1}{2}}(Z_n-Z)$.  In addition, the
increment $X_{n,j}:=M_{n,j}-M_{n,j-1}$ is equal to zero for $j=1$ and,
for $j\geq 2$, it coincides with a random variable of the form
$n^{\gamma-\frac12}(Z_k-Z_{k+1})$ with $k\geq n$. Therefore, we have

\begin{equation*}
\begin{split}
\sum_{j\geq 1} X_{n,j}^2&=n^{2\gamma-1}\sum_{k\geq n}(Z_k-Z_{k+1})^2
=n^{2\gamma-1}\sum_{k\geq n}r_k^2\left(N^{-1}\sum_{i=1}^NI_{k+1}(i)-Z_k\right)^2\\
&\stackrel{a.s.}\sim
c^2 n^{2\gamma-1}\sum_{k\geq n}k^{-2\gamma}
\left(N^{-1}\sum_{i=1}^NI_{k+1}(i)-Z_k\right)^2\\
&\stackrel{a.s.}\longrightarrow
\frac{c^2}{N(2\gamma-1)}Z(1-Z)
\end{split}
\end{equation*}

where the last convergence follows from Lemma~4.1 in
\cite{cri-dai-min} (where, using the notation of such lemma,
$a_k=k^{1-(2\gamma-1)}$,
$b_k=k^{2\gamma-1}$,
$Y_k=\left(N^{-1}\sum_{i=1}^NI_{k+1}(i)-Z_k\right)^2$ and ${\mathcal
  G}_k={\mathcal F}_{k+1}$) and the fact that, by \eqref{forlucca}, we have
\begin{equation}\label{conti}
\begin{split}
\mathrm{E}\left[\left(N^{-1}\sum_{i=1}^N I_{n+1}(i)-Z_n\right)^2 \Big|
\F_n\right]&=
\mathrm{Var}\left[N^{-1}\sum_{i=1}^N I_{n+1}(i) \Big| \F_n\right]\\
&=N^{-2}\sum_{i=1}^N \mathrm{Var}\left[I_{n+1}(i)|\F_n\right]\\
&\stackrel{a.s.}\longrightarrow
N^{-1} Z(1-Z).
\end{split}
\end{equation}

Moreover, again by \eqref{eq-base-inter} (with $\rho=1$), we have

\begin{equation}
\begin{split}
X_n^*=\sup_{j\geq 1}\;|X_{n,j}|&=
n^{\gamma-\frac{1}{2}}\,\sup_{k\geq n}|Z_k-Z_{k+1}|
\leq \sup_{k\geq n}k^{\gamma-\frac{1}{2}}\,|Z_k-Z_{k+1}| \\
&\leq \sup_{k\geq n}k^{\gamma-\frac{1}{2}}\,r_k
\stackrel{a.s.}\longrightarrow 0.
\end{split}
\end{equation}

 Hence, if in Theorem \ref{fam_tri_vet_as_inf} we take $k_n=1$ for
 each $n$ and ${\cal U}$ equal to the $\sigma$-field $\bigvee_n{\cal
   F}_n$, then the conditioning system $({\cal F}_{n,k_n})_{n}$
 coincides with the filtration $\cal F$ and the assumptions are
 satisfied. The proof is thus complete.
\epr

We are now ready for the proof of Theorem \ref{th:fluctuationsZ}.  We
split it into two steps.

\bigskip

\noindent{\em Step 1:} The fact that $\P(Z \in \{0,1\}) < 1$ follows
from Theorem \ref{prop:conv1}. The proof that $\P(Z = z) = 0$ for all
$z \in (0,1)$ is now a consequence of the almost sure conditional
convergence in Proposition \ref{th:stable}, exactly as in Theorem~3.2
in \cite{cri-dai-min}. Indeed, if we denote by $K_n$ a version of the
conditional distribution of $n^{\gamma-\frac12}(Z_n-Z)$ given ${\mathcal F}_n$,
then there exists an event $A$ such that $P(A)=1$ and, for each
$\omega\in A$,
\begin{equation*}
{\textstyle\lim_n} Z_n(\omega)=Z(\omega)\quad\hbox{and}\quad
K_n(\omega)\stackrel{weakly}\longrightarrow
{\mathcal N}\left(0,\frac{c^2}{N(2\gamma-1)}(Z(\omega)-Z^2(\omega))\right).
\end{equation*}

Assume now, by contradiction, that there exists $z\in (0,1)$ with $P(Z=z)>0$,
and set $A'=A\cap\{Z=z\}$ and define $B_n$ as the ${\mathcal
  F}_n$-measurable random set $\{n^{\gamma-\frac12}(Z_n-z)\}$. Then $P(A')>0$
and, since $E\left[I_{\{Z=z\}}\,|\,{\mathcal F}_n\right]$ converges
almost surely to $I_{\{Z=z\}}$, there exists an event $A''$ such that
$P(A'')>0$, $A''\underline{\subset} A'$ and, for each $\omega\in A''$,
\begin{equation*}
K_n(\omega)(B_n(\omega))
\!=\!\E\left[I_{\{n^{\gamma-\frac12}(Z_n-z)\}}\left(n^{\gamma-\frac12}(Z_n-Z)\right)
\,\big|\,{\mathcal F}_n\right](\omega)
\!=\!\E\left[I_{\{Z=z\}}\,|\,{\mathcal F}_n\right](\omega)
\!\longrightarrow\!
I_{\{Z=z\}}(\omega)\!=\!1.
\end{equation*}
On the other hand, we observe that $Z(\omega)-Z^2(\omega)\neq 0$
when $\omega\in A'$. Hence, if $D$ is the discrepancy metric
defined by
$$
D[\mu,\nu]=
{\textstyle\sup_{\{B\in\{\hbox{closed balls of }\R\}\}}} |\mu(B)-\nu(B)|,
$$
which metrizes the weak convergence of a sequence of probability
distributions on $\R$ in the case when the limit distribution is
absolutely continuous with respect to the Lebesgue measure on $\R$
(see \cite{gibbs-2002}), then, for each $\omega\in A'$,  we have
\begin{equation*}
\begin{split}
K_n(\omega)(B_n(\omega))&=
\left|
K_n(\omega)(B_n(\omega))-
{\mathcal N}\left(0,\frac{c^2}{N(2\gamma-1)}(Z(\omega)-Z^2(\omega))
\right)(B_n(\omega))
\right|
\\
&
\leq
D\!\left[
K_n(\omega),
{\mathcal N}\left(0, \frac{c^2}{N(2\gamma-1)}(Z(\omega)-Z^2(\omega))\right)
\right]
\longrightarrow 0.
\end{split}
\end{equation*}
This contradicts the previous fact and the proof of the first step is
concluded.  \\

\noindent{\em Step 2:} We now prove the functional fluctuation result
in \eqref{fcltZ}. First of all, we want to verify the three conditions
$(a2),\, (b2),\, (c2)$ in Theorem \ref{gouet-th} for the stochastic
processes
$$
S^{(n)}_t=
\frac{1}{n^{(2\gamma-1)/2}}
\sum_{k=1}^{\lfloor n t^{1/(2\gamma-1)}\rfloor}
k^{2\gamma-1}(Z_k-Z_{k-1})
$$ in order to obtain a convergence result for $S^{(n)}$ on the space
$(T, m)$ defined in the appendix. Finally, by applying a suitable
continuous transformation, we will arrive to \eqref{fcltZ}.  \\

\noindent{\it Proof of condition (a2):} We want to use Theorem
\ref{dur-res-th}. Let us set ${\mathcal F}_{n,k}={\mathcal F}_k$,
$$
X_{n,k}=\frac{k^{2\gamma-1}(Z_k-Z_{k-1})}{n^{(2\gamma-1)/2}}
\quad\hbox{and}\quad k_n(t)=\lfloor n\, t^{1/(2\gamma-1)} \rfloor
$$
so that $S^{(n)}_t=\sum_{k=1}^{k_n(t)} X_{n,k}$. We observe that

\begin{equation} \label{varcount}
\begin{split}
\sum_{k=1}^{k_n(t)}\mathrm{E}[X_{n,k}^2\,|\,{\mathcal F}_{n,k-1}]&=
\frac{1}{n^{2\gamma -1}}
\sum_{k=1}^{k_n(t)} k^{4\gamma-2}\mathrm{E}[(Z_k-Z_{k-1})^2|{\mathcal F}_{k-1}]\\
&\stackrel{a.s.}\sim
\left(\frac{k_n(t)}{n}\right)^{2\gamma-1}
\frac{c^2}{\left(k_n(t)\right)^{2\gamma-1}}
\sum_{k=1}^{k_n(t)} \frac{1}{k^{1-(2\gamma-1)}}
\mathrm{E}\left[\left(N^{-1}\sum_{i=1}^N I_k(i)-Z_{k-1}\right)^2
\Big|{\mathcal F}_{k-1}\right]\\
&\stackrel{a.s.}\longrightarrow
{\widetilde V}_t= \frac{c^2}{N(2\gamma-1)} Z(1-Z)t
\end{split}
\end{equation}
(where, in the last step, we have used \eqref{conti} and $\lim_n
k_n(t)/n=t^{1/(2\gamma-1)}$) and so condition (a1) of Theorem
\ref{dur-res-th} is verified.  Furthermore, for any $u>2$, we have

\begin{equation*}
\begin{split}
\sum_{k=1}^{k_n(1)}\mathrm{E}[\,|X_{n,k}|^u\,]
&=
\frac{1}{n^{u(2\gamma-1)/2}}\sum_{k=1}^{n}k^{u(2\gamma-1)}
\mathrm{E}[\,|Z_k-Z_{k-1}|^u\,]
=
\frac{1}{n^{u(2\gamma-1)/2}}\sum_{k=1}^{n}k^{u(2\gamma-1)}O(k^{-u\gamma})
\\
&=
\frac{1}{n^{u(2\gamma-1)/2}}\sum_{k=1}^{n}O\left(1/k^{1-[1-u(1-\gamma)]}\right)
=
\begin{cases}
O(n^{-(u/2-1)})\; &\mbox{if } 1-u(1-\gamma)>0\\
O(n^{-u(2\gamma-1)/2})\; &\mbox{if } 1-u(1-\gamma)<0\\
O\!\left(\frac{\ln(n)}{n^{(2\gamma-1)/2(1-\gamma)}}\right)
\; &\mbox{if } 1-u(1-\gamma)=0.
\end{cases}
\end{split}
\end{equation*}
Hence also condition (b1) of Theorem \ref{dur-res-th} holds true (by 
Remark \ref{remark-cond} with $u>2$) and we conclude that
$$
S^{(n)}_t=\sum_{k=1}^{k_n(t)} X_{n,k}=
\frac{1}{n^{(2\gamma-1)/2}}
\sum_{k=1}^{\lfloor n t^{1/(2\gamma-1)}\rfloor}
k^{2\gamma-1}(Z_k-Z_{k-1})
\stackrel{d}\longrightarrow
{\widetilde W}=\left(W_{{\widetilde V}_t}\right)_{t\geq 0}
$$
(w.r.t. Skorohod's topology).
\\

\noindent{\it Proof of condition (b2):}
For each $\epsilon>0$,  we observe that we have

\begin{equation*}
\begin{split}
\mathrm{E}[ n^{ 2(1/2+\epsilon)(1-2\gamma) }
n^{4\gamma-2} (Z_n-Z_{n-1})^2]
\leq
n^{ 2(1/2+\epsilon)(1-2\gamma)+ 4\gamma-2}r_{n-1}^2
&\sim
c^2
n^{ 2(1/2+\epsilon)(1-2\gamma)+ 4\gamma-2 -2\gamma}
\\
&=
\frac{c^2}{ n^{1+ 2\epsilon(2\gamma-1) } }.
\end{split}
\end{equation*}
Therefore the martingale

\begin{equation*}
\left(\sum_{k=1}^{n} k^{(1/2+\epsilon)(1-2\gamma)} k^{(2\gamma-1)} (Z_k-Z_{k-1})\right)
\end{equation*}
is bounded in $L^2$ and so $\sum_{k=1}^{+\infty}
k^{(1/2+\epsilon)(1-2\gamma)} k^{(2\gamma-1)} (Z_k-Z_{k-1})$
is a.s. convergent. By Kronecker's lemma, we get $\sum_{k=1}^{n}
k^{(2\gamma-1)}(Z_k-Z_{k-1})=o(n^{(1/2+\epsilon)(2\gamma-1)})$ a.s.
This fact implies that, for each fixed $n$, the process $S^{(n)}$ is
such that $S^{(n)}_t=o(t^{(1/2+\epsilon)})$ a.s. as $t\to +\infty$.
\\

\noindent{\it Proof of condition (c2):} Fix $\theta>1/2$, $\epsilon>0$
and $\eta>0$. We want to verify that there exists $t_0$ such that

\begin{equation}\label{inequality}
P\left\{
\sup_{t\geq t_0}\frac{ |S^{(n)}_t| }{t^{\theta}}
>\epsilon
\right\}
\leq\eta.
\end{equation}
To this purpose, we observe that

$$
S^{(n)}_t=\frac{1}{n^{(2\gamma-1)/2}}L_{\lfloor nt^{1/(2\gamma-1)}\rfloor}\quad
\hbox{where }
L_k=\sum_{j=1}^k j^{2\gamma-1}(Z_j-Z_{j-1})=
\sum_{j=1}^{k} \xi_j.
$$
Denoting by $C$ a suitable positive constant (which may vary at each step),
we have

\begin{equation*}
\begin{split}
P\left\{\sup_{t\geq (n_0/n)^{2\gamma-1}}\frac{ |S^{(n)}_t| }{ t^{\theta}}>
\epsilon\right\}
&\leq
P\left\{|L_k|>
\frac{\epsilon k^{(2\gamma-1)\theta}}{2n^{(2\gamma-1)\theta-(2\gamma-1)/2}}
\;\hbox{for some }k\geq n_0\right\}
\\
&\leq
\sum_{i=1}^{+\infty}
P\left\{\max_{2^{i-1}n_0\leq k\leq 2^i n_0}|L_k|>
\frac{\epsilon \left(2^{i-1}n_0\right)^{(2\gamma-1)\theta}}
{2n^{(2\gamma-1)\theta-(2\gamma-1)/2}}\right\}
\\
&\leq
\frac{16 n^{2(2\gamma-1)\theta-(2\gamma-1)}}
{\epsilon^2 n_0^{2(2\gamma-1)\theta}}
\sum_{i=1}^{+\infty} 2^{-2i}
\sum_{j=1}^{2^in_0} {\mathrm{E}}[\xi_j^2]
\\
&\leq
C\frac{n^{2(2\gamma-1)\theta-(2\gamma-1)}}
{\epsilon^2 n_0^{2(2\gamma-1)\theta}}
\sum_{i=1}^{+\infty} 2^{-2i}(2^{i}n_0)^{2\gamma-1}\\
&=
\frac{C}{\epsilon^2}\left(\frac{n}{n_0}\right)^{(2\gamma-1)(2\theta-1)}
\sum_{i=1}^{+\infty}\left(\frac{1}{2^{3-2\gamma}}\right)^i=
\frac{C}{\epsilon^2}\left(\frac{n}{n_0}\right)^{(2\gamma-1)(2\theta-1)}
\end{split}
\end{equation*}
where the third inequality is the H\'ajek-R\'enyi inequality for
martingales (see the appendix) and we used the fact that
$\mathrm{E}[\xi_j^2]\sim C/j^{1-(2\gamma-1)}$ with $C>0$.  Therefore,
in order to obtain (\ref{inequality}), it is enough to set
$t_0^{2\theta-1}=\frac{C}{\epsilon^2\eta}$.  \\

\noindent{\it Conclusion:} By Theorem \ref{gouet-th}, we have
$$
S^{(n)}\stackrel{d}\longrightarrow \widetilde{W}=
\left(W_{\widetilde{V}_t}\right)_{t\geq 0}
\qquad\hbox{on } (T, m).
$$
Now, let $g: T\to T_1^*$ be the Barbour's transform
(defined in the appendix) and observe that, since

\begin{equation*}
(\Delta S^{(n)})_s=S^{(n)}_s-S^{(n)}_{s-}=
\begin{cases}
\frac{k^{2\gamma-1}(Z_k-Z_{k-1})}{n^{(2\gamma-1)/2}}
\quad &\hbox{if } k=ns^{1/(2\gamma-1)}\\
0 \quad &\hbox{otherwise},
\end{cases}
\end{equation*}

we have
\begin{equation*}
\begin{split}
g(S^{(n)})(t)&=
\sum_{s\geq 1/t}\frac{(\Delta S^{(n)})_s}{s}
=\sum_{k> \lfloor n t^{-1/(2\gamma-1)}\rfloor}
\left(\frac{n}{k}\right)^{2\gamma-1}
\frac{k^{2\gamma-1}(Z_k-Z_{k-1})}{n^{(2\gamma-1)/2}}
\\
&=
n^{(2\gamma-1)/2}
\sum_{k> \lfloor nt^{-1/(2\gamma-1)}\rfloor}
(Z_k-Z_{k-1})
=
n^{(2\gamma-1)/2}
\left(Z-Z_{\lfloor n t^{-1/(2\gamma-1)}\rfloor}\right).
\end{split}
\end{equation*}

Therefore, by the properties of $g$ (see the appendix),  we get
$$
g(S^{(n)})=\left(n^{(2\gamma-1)/2}
\left(Z-Z_{\lfloor n t^{-1/(2\gamma-1)}\rfloor}\right)\right)_{t\geq 0}
\stackrel{d}\longrightarrow g(\widetilde W)\stackrel{d}=\widetilde W
\qquad\hbox{on } (T_1^*, m_1^*).
$$ Immediately, by symmetry of $W$, we get the convergence result of
  $\left(n^{(2\gamma-1)/2}(Z_{\lfloor n
    t^{-1/(2\gamma-1)}\rfloor}-Z)\right)_{t\geq 0}$ to the stochastic
  process $\widetilde{W}$. Then, we can apply the continuous map $h:
  T_1^*\to D$ defined as
\begin{equation*}
h(f)(0)=0\quad\mbox{and}\quad h(f)(t)=tf(t^{-1})
\end{equation*}
(see \cite{hida}) and obtain

\be{fcltZ-intermedio}
\left(t n^{\g- \frac12} (Z_{\lfloor n t^{1/(2\gamma-1)}\rfloor}-Z)
\right)_{t\geq 0}
\stackrel{d}\longrightarrow
\left(h({\widetilde W})\right)_{t\geq 0}
\stackrel{d}= {\widetilde W}=
\left(W_{{\widetilde V}_t}\right)_{t\geq 0}
\qquad\hbox{(w.r.t. Skorohod's topology)}.
\ee

Finally, we can set $t=s^{2\gamma-1}$ and obtain
\begin{equation*}
\left(s^{2\gamma-1} n^{\g- \frac12}(Z_{\lfloor n s \rfloor}-Z)
\right)_{s\geq 0}
\stackrel{d}\longrightarrow
\left(W_{V_s}\right)_{s\geq 0}
\qquad\hbox{(w.r.t. Skorohod's topology)},
\end{equation*}
which coincides with \eqref{fcltZ}.

\section{Fluctuation theorems II: proof of Theorem \ref{FCLTSyncro}}

In all the sequel, we denote by $C>0$ a suitable constant (which may
vary at each step).
\\

By \eqref{eq-base-diff}, recalling that we are
assuming $\rho = 1$, we have
\begin{equation*}
\E[Z_{n+1}(i)-Z_{n+1}|\F_n] = (1-\alpha  r_n)\left(Z_n(i) - Z_n\right).
\end{equation*}
Thus, setting
\begin{equation*}
l_n = \prod_{k=0}^{n-1} \frac{1}{1-\alpha  r_k},
\end{equation*}
so that $l_n = (1-\alpha  r_n) l_{n+1}$, we have that $L_n = l_n
\left(Z_n(i)- Z_n\right)$ forms a martingale with
\begin{equation*}
\xi_n = \Delta L_n= L_n - L_{n-1} =
l_n \left( Z_n(i)-Z_n - \E[Z_n(i)-Z_n|\F_{n-1}] \right).
\end{equation*}
Observe that, if we fix $\epsilon\in (0,1)$, for all $n\geq \bar{n}$
and $\bar{n}$ sufficiently large
\begin{equation*}
l_{\bar{n}}\exp\left[-\sum_{k=\bar{n}}^{n}\ln\left(1-\frac{\alpha
    c_1}{k^{\g}} \right)\right] \leq l_n\leq
l_{\bar{n}}\exp\left[-\sum_{k=\bar{n}}^{n}\ln\left(1-\frac{\alpha
    c_2}{k^{\g}} \right)\right]
\end{equation*}
where $c_1=c(1-\epsilon)$ and $c_2=c(1+\epsilon)$.

In particular, recalling that the function $\varphi(x)= - x -
\ln(1-x)$ is such that $0 \leq \varphi(x) \leq C x^2$ for $0 \leq x
\leq \bar{x}<1$, and assuming $\bar n$ large enough such that $\alpha
c_2/k^\g \leq \bar{x}$ for $k\geq \bar{n}$, we have for $n\geq\bar{n}$

\begin{equation}\label{disuguaglianza1}
l_{\bar{n}}\exp\left[\sum_{k=\bar{n}}^{n}\frac{\alpha  c_1}{k^\g}\right]
\leq l_n\leq l_{\bar{n}}\exp\left[\sum_{k=\bar{n}}^{n}\left(\frac{\alpha  c_2}{k^\g}
+C \frac{(\alpha  c_2)^2}{k^{2\g}}\right)\right].
\end{equation}
As a simple consequence, we have for $n$ large enough
\begin{equation}\label{lnlim1}
l_n \geq l_{\bar{n}}\exp\left[\sum_{k=\bar{n}}^{n}\frac{\alpha
    c_1}{k^\g}\right]\geq \left\{
\begin{array}{ll}
C
\exp \left( \frac{c_1\alpha }{1-\g} n^{1-\g} \right) & \mbox{for } 1/2<\g<1 \\ 
C n^{c_1 \alpha} & \mbox{for } \g=1.
\end{array}
\right.
\end{equation}
This, in particular, implies
\be{suminf}
\lim_n\,\frac{n^{\g}}{l_n^{2}}= 0\qquad\mbox{and}\qquad 
\sum_{n=1}^{+\infty} l_n^2 r_{n-1}^2 = +\infty.  
\ee 
Indeed, these facts follow from \eqref{lnlim1} immediately for
$1/2<\g<1$; while, for $\g =1$ one has to note that, since we assume
$2c\alpha >1$, we can choose $\e$ small enough so that $c_1 \alpha >
\frac12$.

We now use Theorem \ref{dur-res-th} to obtain a functional central
limit theorem for $(L_n)$, from which the corresponding result for
$(Z_n(i)-Z_n)$ will follow. Set

\begin{equation*}
X_{n,k} = \frac{n^{\g/2} \xi_k}{l_n},
\end{equation*}
$\F_{n,k} = \F_k$ and $k_n(t) = \lfloor n + n^{\g} t \rfloor$.  We
start with showing condition (a1) of Theorem \ref{dur-res-th}. Note that

\begin{equation}\label{conda1}
\begin{split}
\sum_{k=1}^{k_n(t)}{\mathrm{E}}[X_{n,k}^2\,|\,{\cal F}_{n,k-1}] & =
\frac{n^{\g}}{l^2_n} \sum_{k=1}^{k_n(t)}\E\left[\xi_k^2|\F_{k-1} \right]
=
\frac{n^{\g}}{l^2_n} \sum_{k=1}^{k_n(t)} l^2_k \mathrm{Var}[Z_k(i)-Z_k|\F_{k-1}]
\\
& \stackrel{a.s.}\sim
(1-1/N)Z(1-Z)
\frac{n^{\g}}{l^2_n} \sum_{k=1}^{k_n(t)} l^2_k r^2_{k-1},
\end{split}
\end{equation}
where, in the last step, we have used the fact that $\lim_n
n^{\g}/l^2_n=0$ and \eqref{sim-var} with $\rho = 1$. In order to 
estimate the sum $\sum_{k=1}^{k_n(t)} l^2_k r^2_{k-1}$, we observe
that, for $k\to +\infty$, we have
\begin{eqnarray}\label{lkrk}
\frac{1}{k^{\g}} l^2_k - \frac{1}{(k-1)^{\g}} l^2_{k-1} &=&
\left[\frac{1}{k^{\g}} - \frac{1}{(k-1)^{\g}} \right] l^2_{k-1}
+ \frac{1}{k^{\g}}\left(l^2_k -  l^2_{k-1}\right)\nonumber \\
&=&
\left[-\frac{\gamma}{k^{\g + 1}}+o(1/k^{\gamma+1})\right] 
l^2_k \frac{l^2_{k-1}}{l^2_k}
+
\frac{1}{k^{\g}} l^2_k \left( 1- \frac{l^2_{k-1}}{l^2_k} \right) \nonumber\\
&=&
-\frac{\gamma}{k^{\g + 1}} l^2_k (1-\alpha r_{k-1})^2  
+ o(l_k^2/k^{\gamma+1})
+ \frac{1}{k^{\g}} l^2_k
\left[ 1 - (1-\alpha  r_{k-1})^2 \right]\nonumber \\
&=&
-\frac{\gamma}{k^{\g + 1}} l^2_k  +
\frac{1}{k^\g} l^2_k
\left( - \alpha^2 r_{k-1}^2 + 2\alpha  r_{k-1} \right)
+o(l_k^2/k^{\gamma+1})
\nonumber
\\
&=& \left\{\begin{array}{ll}
\left(\frac{2\alpha  r_{k-1}}{k}-\frac{1}{k^2}\right) l^2_k + o(l_k^2 r_{k-1}^2)
&\mbox{if } \g=1,\, 2c\alpha >1\\
\frac{2\alpha  r_{k-1}}{k^{\g}} l^2_k + o(l_k^2 r_{k-1}^2) &\mbox{if } 1/2<\g<1\\
\end{array}\right.\nonumber\\
&\sim& \left\{\begin{array}{ll}
\frac{2c \alpha -1}{c^2}l_k^2 r_{k-1}^2  &\mbox{if } \g=1,\, 2c\alpha >1\\
\frac{2\alpha }{c} l^2_k r_{k-1}^2  &\mbox{if } 1/2<\g<1.\\
\end{array}\right.
\end{eqnarray}

Therefore, recalling \eqref{suminf},  we have
\begin{equation}\label{est0}
\frac{n^{\g}}{l^2_n}
\sum_{k=1}^{k_n(t)} l^2_k r^2_{k-1}
\sim
\begin{cases}
\frac{c^2}{2c\alpha-1 }\frac{n}{l^2_n} \sum_{k=2}^{k_n(t)}
\left( \frac{1}{k} l^2_k - \frac{1}{(k-1)} l^2_{k-1}\right)
\sim \frac{c^2}{2c\alpha-1 }
\frac{n}{k_n(t)} \frac{l^2_{k_n(t)}}{l^2_n}\;&\mbox{if } \gamma=1,\, 2c\alpha >1 
\\
\frac{c}{2\alpha }\frac{n^{\g}}{l^2_n} \sum_{k=2}^{k_n(t)}
\left( \frac{1}{k^{\g}} l^2_k - \frac{1}{(k-1)^{\g}} l^2_{k-1}\right)
\sim \frac{c}{2\alpha }
\frac{n^{\g}}{k_n(t)^{\g}} \frac{l^2_{k_n(t)}}{l^2_n}\;&\mbox{if } 1/2<\gamma<1.
\end{cases}
\end{equation}

Observing that
\begin{equation}\label{uno}
\lim_{n} \left(\frac{n}{k_n(t)}\right)^\g = \left\{ \begin{array}{ll}
  1 & \mbox{if } 1/2<\g<1 \\ \frac{1}{1+t} & \mbox{if } \g=1, \end{array}
\right.
\end{equation}
we are left to compute the limit
\begin{equation*}
\lim_{n}  \frac{l_{k_n(t)}}{l_n}.
\end{equation*}
This is done in the following Lemma.

\bl{lemma:timechange} The following convergence holds uniformly over
compact subsets of $[0,+\infty)$:
\begin{equation}\label{timechange}
\lim_{n} \frac{l_{k_n(t)}}{l_n}= \left\{ \begin{array}{ll} e^{c\alpha
    t} & \mbox{for } 1/2<\g<1 \\ 
(1+t)^{c\alpha} & \mbox{for } \g=1,\, 2c\alpha >1. 
\end{array}
\right.
\end{equation}
\el

We postpone the proof of Lemma \ref{lemma:timechange} and complete
the proof of Theorem \ref{FCLTSyncro}. Inserting \eqref{uno} and
\eqref{timechange} in \eqref{est0} and \eqref{conda1} we obtain

\begin{equation*}
\sum_{k=1}^{k_n(t)}{\mathrm{E}}[X_{n,k}^2\,|\,{\cal F}_{n,k-1}]
\stackrel{a.s.}\longrightarrow \left\{ \begin{array}{ll}
  \left(1-1/N\right)\frac{c}{2\alpha } e^{2c\alpha  t} Z(1-Z) & \mbox{if }
  1/2<\g<1 \\ \left(1-1/N\right)\frac{c^2}{2c\alpha -1}(1+t)^{2c\alpha-1} Z(1-Z)
  & \mbox{if } \g=1,\, 2c\alpha >1. \end{array} \right.
\end{equation*}
and condition (a1) of Theorem \ref{dur-res-th} is verified.\\

We now prove condition (b1) of Theorem \ref{dur-res-th}, via Remark
\ref{remark-cond} (with $u>2$) and the sufficient condition in
\eqref{suff}. Let $u>2$.  Note that, using \eqref{eq-base-diff} and
\eqref{s_cond_zni-zn} with $\rho=1$, we have $|\xi_{n+1}|\leq 2
l_{n+1} r_n$ for all $n$.  Then
\begin{equation}\label{est_u}
\sum_{k=1}^{k_n(1)}
{\mathrm{E}}[\,|X_{n,k}|^u\,]  \leq
2^u \frac{n^{\frac{\g u}{2}}}{l_n^u} \sum_{k=1}^{k_n(1)}  l_k^u r_{k-1}^u.
\end{equation}
This last sum can be estimated as we have done in \eqref{est0}.
A computation analogous to the one made in \eqref{lkrk} gives
\begin{eqnarray}\label{lkrk_u}
\frac{1}{k^{\g(u-1)}} l^u_k - \frac{1}{(k-1)^{\g(u-1)}} l^u_{k-1}
&\sim& \left\{\begin{array}{ll}
\frac{uc \alpha -(u-1)}{c^u}l_k^u r_{k-1}^u  &\mbox{if } \g=1,\, uc\alpha>(u-1)
\\
\frac{u\alpha }{c^{u-1}} l^u_k r_{k-1}^u  &\mbox{if } 1/2<\g<1.\\
\end{array}\right.
\end{eqnarray}
Then, for $1/2<\gamma\leq 1$ and suitable $u>2$ (note that, when
$\g=1$, we have to choose $u>2$ such that
$uc\alpha-(u-1)=u(c\alpha-1)+1>0$ and this choice is possible since in this
case $2c\alpha>1$ by assumption), we have
$$
\frac{n^{\frac{\g u}{2}}}{l_n^u} \sum_{k=1}^{k_n(1)}  l_k^u r_{k-1}^u \sim
C\left(\frac{l_{k_n(1)}}{l_n}\right)^u
\frac{n^{\frac{\g u}{2}}}{n^{\g(u-1)}}\longrightarrow 0.
$$
Thus, by applying Theorem \ref{dur-res-th}, we get
\begin{equation*}
\left(\frac{n^{\g/2}}{l_n} L_{k_n(t)} \right)_{t \geq 0}
\stackrel{d}\longrightarrow \left( W_{V_t} \right)_{t \geq 0},
\end{equation*}
where $(V_t)_{t\geq 0}$ and $(W_t)_{t \geq 0}$ are defined in the
statement of the theorem. Finally, recalling that
$\frac{n^{\g/2}}{l_n} L_{k_n(t)}=n^{\g/2}\frac{l_{k_n(t)}}{l_n}
\left(Z_{k_n(t)}(i) - Z_{k_n(t)} \right)$ and using Lemma \ref{lemma:timechange}
the proof is complete.
\qed

\bigskip

\noindent
{\em Proof of Lemma \ref{lemma:timechange}}. We begin by proving the
claim for $1/2<\g<1$. Note first that $\frac{l_{k_n(t)}}{l_n}=1$ if $t<
n^{-\gamma}$, then for any $T\geq 1$ we have

\begin{equation}
\sup_{t\in [0,T]}\left | \frac{l_{k_n(t)}}{l_n}-e^{c\alpha  t}\right |=
\max \left\{ \sup_{t\in [0,n^{-\gamma})}|1-e^{c\alpha  t}|,
\sup_{t\in [n^{-\gamma},T]}e^{c\alpha  t}
\left | \frac{l_{k_n(t)}}{l_n}e^{-c\alpha  t}-1\right | \right\}.
\end{equation}
Since $ \sup_{t\in [0,n^{-\gamma})}|1-e^{c\alpha  t}|=e^{c\alpha
    n^{-\gamma}}-1$ converges to 0, it is enough to show that
$$
\lim_n\, \sup_{t\in [n^{-\gamma},T]}
\left | \frac{l_{k_n(t)}}{l_n}e^{-c\alpha t}-1\right |=0.
$$

We fix $\sigma>0$ and take
$0<\epsilon<\min\{\frac{\ln(1+\sigma)}{2c\alpha T}, 1\}$. Then, using
\eqref{disuguaglianza1}, we have, for all $n\geq \bar{n}$ and
$\bar{n}$ sufficiently large,

\begin{equation}\label{disuguaglianza}
\exp\left[\sum_{k=n}^{k_n(t)}
\frac{\alpha  c_1}{k^{\g}}\right]
 \leq
\frac{l_{k_n(t)}}{l_n}\leq
\exp\left[\sum_{k=n}^{k_n(t)}\left(\frac{\alpha  c_2}{k^{\g}}
+C\frac{(\alpha  c_2)^2}{k^{2\g}}\right)\right]\ \ \ \ \forall t\in [n^{-\g},T]
\end{equation}
where $c_1=c(1-\epsilon)$ and $c_2=c(1+\epsilon)$.

Now we use the following asymptotics, which
hold for all $p>0$ and $n,m\geq 1$:
\begin{equation}\label{asy-somma}
\sum_{k=n}^{n+m-1} \frac{1}{k^p} =
\left\{ \begin{array}{ll}
\e_p(n,m) & \mbox{if } p>1 \\
\ln \left(1+\frac{m}{n}\right) + \e_1(n,m)
& \mbox{if } p=1 \\
\frac{n^{1-p}}{1-p} \left[ \left(1+\frac{m}{n} \right)^{1-p} - 1 \right]
+ \e_p(n,m) & \mbox{if } p<1
\end{array} \right.
\end{equation}
where $\e_p(n,m)$ denotes a positive function such that $\lim_{n}
\sup_{m \geq 1} \e_p(n,m) =0$. Thus, recalling that $1/2<\g<1$,

\begin{equation}\label{dis1}
\exp\left\{\frac{\alpha  c_1 n^{1-\g}}{1-\g} \left[\left(1+\frac{\lfloor
    n^\g t\rfloor}{n} \right)^{1-\g} - 1 \right]\right\} \leq
\frac{l_{k_n(t)}}{l_n} \leq \exp\left\{\frac{\alpha  c_2 n^{1-\g}}{1-\g}
\left[\left(1+\frac{\lfloor n^\g t\rfloor}{n} \right)^{1-\g} - 1
  \right]+ A_{n,t}\right\}
\end{equation}
where $A_{n,t}=\alpha  c_2\e_\g(n,\lfloor n^\g t\rfloor)+C(\alpha  c_2)^2
\e_{2\g}(n,\lfloor n^\g t\rfloor)$.  Moreover, using the relation
$1+(1-\gamma)x-\frac{(1-\gamma)\gamma}{2}x^2\leq (1+x)^{1-\gamma} \leq
1+ (1-\gamma)x$ for all $x\geq 0$ and $1/2<\gamma<1$, we have that
\begin{equation}\label{dis2}
t-\frac{1}{n^{\gamma}}-\frac{\gamma t^2}{2 n^{1-\gamma}}\leq
 \frac{n^{1-\g}}{1-\g}
\left[ \left(1+\frac{\lfloor n^\g t\rfloor}{n} \right)^{1-\g} - 1 \right]
\leq t
\end{equation}
 Therefore, using \eqref{dis1} and \eqref{dis2} we obtain
\[
e^{-c\alpha  \epsilon t 
- \frac{c\alpha (1-\epsilon)}{n^\gamma} - \frac{c\alpha (1-\epsilon)\gamma t^2}{2
    n^{1-\gamma}}}\leq \frac{l_{k_n(t)}}{l_n}e^{-c\alpha  t}\leq
e^{c\alpha \epsilon t+ A_{n,t}} \ \ \ \ \forall t\in [n^{-\g},T]
\]
and, since $\lim_n \sup_t A_{n,t}=0$, for all $n$ sufficiently large we have
\[
e^{-2c\alpha  \epsilon T}\leq \frac{l_{k_n(t)}}{l_n}e^{-c\alpha  t}
\leq e^{2c\alpha \epsilon T}
\ \ \ \ \forall t\in [n^{\g},T]
\]
from which it follows
\begin{equation*}
\sup_{t\in [n^{-\gamma},T]}\left | \frac{l_{k_n(t)}}{l_n}e^{-c\alpha  t}-1\right |\leq
e^{2c\alpha \epsilon T}-1<\sigma
\end{equation*}
for all $n$ sufficiently large and the proof is complete, since
$\sigma$ can be taken arbitrarily small.  \\

\medskip
We now turn to the case $\g=1$ and $2c\alpha >1$. As for $1/2<\g<1$,
it is enough to show that $\sup_{t\in [n^{-1},T]}\left |
\frac{l_{k_n(t)}}{l_n}(1+t)^{-c\alpha }-1\right | $ converges to
0.\ \ We fix $\sigma>0$ and take
$0<\epsilon<\min\{\frac{\ln(1+\sigma)}{2c\alpha\ln(1+T)},
1\}$. Using \eqref{disuguaglianza} and \eqref{asy-somma} for $\g=1$
, we obtain, for all $n\geq \bar{n}$ and $\bar{n}$ sufficiently large,
\begin{equation}\label{dis1_2}
\exp\left[\alpha  c_1 \ln \left(1+\frac{\lfloor n t\rfloor
  }{n}\right)\right] \leq \frac{l_{k_n(t)}}{l_n} \leq \exp\left[\alpha
  c_2 \ln \left(1+\frac{\lfloor n t\rfloor }{n}\right) +
  A_{n,t}\right]\ \ \forall t \in [n^{-1}, T]
\end{equation}
where $c_1=c(1-\epsilon)$,\ $c_2=c(1+\epsilon)$ and $A_{n,t}=\alpha
c_2\e_1(n,\lfloor n t\rfloor)+C(\alpha  c_2)^2 \e_{2}(n,\lfloor n
t\rfloor)$.  Then, for all $n$ sufficiently large
\[ \left(1+\frac{\lfloor nt \rfloor}{n}\right)^{c(1-\epsilon)\alpha }(1+t)^{-c\alpha }
 \leq \frac{l_{k_n(t)}}{l_n}(1+t)^{-c\alpha }\leq \left(1+\frac{\lfloor
   nt
   \rfloor}{n}\right)^{c(1+\epsilon)\alpha }e^{A_{n,t}}(1+t)^{-c\alpha }
 \ \ \ \ \forall t\in [n^{-1},T]
\]
and, using the fact that $\lim_n \sup_t A_{n,t}=0$, for all $n$
sufficiently large we can write
\[
(1+T)^{-2c\alpha \epsilon}\leq \frac{l_{k_n(t)}}{l_n}(1+t)^{-c\alpha }\leq
(1+T)^{2c\alpha \epsilon} \ \ \ \ \forall t\in [n^{-1},T]
\]
from which it follows
\begin{equation*}
\sup_{t\in [n^{-1},T]}\left | \frac{l_{k_n(t)}}{l_n}(1+t)^{-c\alpha }-1\right |\leq
(T+1)^{2c\alpha \epsilon}-1<\sigma
\end{equation*}
for all $n$ sufficiently large and the
proof is complete, since $\sigma$ can be taken arbitrarily small.
\qed

\section{Fluctuation theorems III: proof of Theorems \ref{FCLT-q} and
\ref{FCLTSyncro-q}}

\subsection{Proof of Theorem \ref{FCLT-q}}

\noindent {\em Proof of case (i)}.  The proof is the same as for
Theorem \ref{FCLTSyncro}, part (i). Indeed,
using \eqref{s_cond_zn-q} we can define
$l_n=\prod_{k=0}^{n-1}\frac{1}{1-(1-\rho)  r_k}$ and
$X_{n,k}=\frac{n^{\g/2}\xi_k}{l_n}$, where
$\xi_n=l_n(Z_n-\mathrm{E}[Z_n|\mathcal{F}_{n-1}])$. Using the estimate
\eqref{suminf} and \eqref{sim-var-zn} we obtain \be{conda1_2_q}
\begin{split}
\sum_{k=1}^{k_n(t)}{\mathrm{E}}[X_{n,k}^2\,|\,{\cal F}_{n,k-1}] &
\stackrel{a.s.}\sim \frac{1}{N}q(1-q)\rho^2\frac{n}{l^2_n}\sum_{k=1}^{k_n(t)}
l^2_k r^2_{k-1}.
\end{split}
\ee Then, using estimate \eqref{est0} and Lemma \ref{lemma:timechange} with
$(1-\rho)$ instead of $\alpha$, we can conclude that
$\frac{n}{l^2_n}\sum_{k=1}^{k_n(t)} l^2_k r^2_{k-1}$ converges to
$\frac{c}{2(1-\rho)} e^{2c(1-\rho)t}$ and so condition (a1) of Theorem
\ref{dur-res-th} is satisfied.\\

Using \eqref{eq-base-inter} and \eqref{s_cond_zn-q} we have
$|\xi_{n+1}|\leq l_{n+1}r_n$ and so, fixing $u>2$, we can repeat the
same argument used in the proof of Theorem \ref{FCLTSyncro} in order
to prove
\begin{equation}\label{est_u_q}
\sum_{k=1}^{k_n(1)}
{\mathrm{E}}[\,|X_{n,k}|^u\,]\longrightarrow 0
\end{equation}
and the conclusion follows from Theorem \ref{dur-res-th}, Remark
\ref{remark-cond} (with $u>2$) and Lemma \ref{lemma:timechange}.\\

\noindent {\em Proof of case (ii).}  Following proof and notations of
Theorem \ref{FCLTSyncro}, part (ii), with $1-\rho$ in place of
$\alpha$, since $2c(1-\rho)>1$, relation \eqref{conda1} becomes
\begin{equation*}
\begin{split}
\sum_{k=1}^{k_n(t)}{\mathrm{E}}[X_{n,k}^2\,|\,{\cal F}_{n,k-1}] &
\stackrel{a.s.}\sim
\frac{1}{N}q(1-q)\rho^2\frac{n}{l^2_n}\sum_{k=1}^{k_n(t)} l^2_k
r^2_{k-1}.
\end{split}
\end{equation*}
and, by estimate \eqref{est0} and Lemma \ref{lemma:timechange} with
$(1-\rho)$ instead of $\alpha$, the last sum converges to
$\frac{c^2}{2c(1-\rho)-1}(1+t)^{2c(1-\rho)-1}$. Analogously,
condition (b1) of Theorem \ref{dur-res-th} and then the conclusion can
be easily derived as seen in Theorem \ref{FCLTSyncro}.

\subsection{Proof of Theorem \ref{FCLTSyncro-q}}

\noindent {\em Proof of case (i)}.  We again follow the proof of
Theorem \ref{FCLTSyncro} by posing $\beta =1-(1-\alpha)\rho$, using
\eqref{s_cond_zni-zn} and defining
$l_n=\prod_{k=0}^{n-1}\frac{1}{1-\beta  r_k}$. Then, given
$X_{n,k}=\frac{n^{\g/2}\xi_k}{l_n}$, where
$\xi_n=l_n(Z_n(i)-Z_n-\mathrm{E}[Z_n(i)-Z_n|\mathcal{F}_{n-1}])$, we
use the estimate \eqref{suminf} and \eqref{sim-var} and we have
\begin{equation*}
\begin{split}
\sum_{k=1}^{k_n(t)}{\mathrm{E}}[X_{n,k}^2\,|\,{\cal F}_{n,k-1}] &
\stackrel{a.s.}\sim (1-\frac{1}{N})q(1-q)\rho^2\frac{n}{l^2_n}\sum_{k=1}^{k_n(t)}
l^2_k r^2_{k-1}.
\end{split}
\end{equation*}
Then, by estimate \eqref{est0} and Lemma \ref{lemma:timechange} with
$\beta$ instead of $\alpha$, it follows
$\frac{n}{l^2_n}\sum_{k=1}^{k_n(t)} l^2_k r^2_{k-1}\longrightarrow
\frac{c}{2(1-(1-\alpha)\rho )} e^{2c(1-(1-\alpha)\rho )t}$ and
condition (a1) of Theorem \ref{dur-res-th} is verified.\\ For
condition (b1), by \eqref{eq-base-diff} and \eqref{s_cond_zni-zn}, we
obtain $|\xi_{n+1}|\leq 2 l_{n+1}r_n$ and so, for a fixed $u>2$
\begin{equation*}
\sum_{k=1}^{k_n(1)}
{\mathrm{E}}[\,|X_{n,k}|^u\,]\longrightarrow 0
\end{equation*}
and the conclusion follows from Theorem \ref{dur-res-th}, Remark
\ref{remark-cond} (with $u>2$) and Lemma \ref{lemma:timechange}.\\

\noindent {\em Proof of case (ii).}  The proof is
analogous to the one of case (ii) of Theorem \ref{FCLT-q} above.
\\

\bigskip

\noindent {\bf Acknowledgment.} \\

\noindent
The authors thank the Gran Sasso Science Institute for
having partially supported their collaboration.

\noindent Irene Crimaldi and Ida G. Minelli are members of the Italian
group ``Gruppo Nazionale per l'Analisi Matematica, la
Pro\-ba\-bi\-li\-t\`a e le loro Applicazioni (GNAMPA)'' of the Italian
Institute ``Istituto Nazionale di Alta Matematica (INdAM)''.

\noindent Irene Crimaldi acknowledges support from CNR PNR Project
``CRISIS Lab''.

\appendix

\section{A technical lemma}

\begin{lemma}\label{lemma-tecnico}
Let $(x_n)$ be a sequence of positive numbers that satisfies the
following equation:
\begin{equation}\label{equation-y}
x_{n+1}=(1-a r_n)x_n + K_n r_n^2
\end{equation}
where $a>0$, $r_n \geq 0$ and $0\leq K_n \leq K$. Suppose that
\be{app-r_n}
\sum_{n} r_n = +\infty \ \ \ \mbox{and} \ \ \
\sum_{n} r_n^2 < +\infty.
\ee
Then
\[
\lim_{n \ra +\infty} x_n = 0.
\]
\end{lemma}

{\bf Proof.}  The case $K=0$ is well-known and so we will prove the
statement with $K>0$. Let $l$ be such that $ar_n < 1$ for all $n \geq
l$. Then for $n \geq l$ we have $x_n \leq y_n$, where
\[
\left\{ \begin{array}{rcl}
y_{n+1} & = &(1-a r_n)y_n + K r_n^2  \\ y_l & = & x_l \end{array} \right.
\]
Set $\varepsilon_n=ar_n$ and $\delta_n=K r_n^2$. It holds
$$ y_n= y_l \prod_{i=l}^{n-1}(1-\varepsilon_i)+ \sum_{i=l}^{n-1}
\delta_i \prod_{j=i+1}^{n-1} (1-\varepsilon_j).
$$
Using the fact that $\sum_n \varepsilon_n = +\infty$, it follows that
\[
\prod_{i=l}^{n-1}(1-\varepsilon_i) \longrightarrow 0.
\]
Moreover, for every $m \geq l$,
\begin{eqnarray}
\sum_{i=l}^{n-1} \delta_i \prod_{j=i+1}^{n-1} (1-\varepsilon_j) & = &
\sum_{i=l}^{m-1} \delta_i \prod_{j=i+1}^{n-1} (1-\varepsilon_j) +
\sum_{i=m}^{n-1} \delta_i \prod_{j=i+1}^{n-1}
(1-\varepsilon_j) \label{app-y_n} \\ & \leq & \prod_{j=m}^{n-1}
(1-\varepsilon_j) \sum_{i=l}^{m-1} \delta_i + \sum_{i=m}^{+\infty}
\delta_i. \nonumber
\end{eqnarray}

Using the fact that $ \prod_{j=m}^{n-1} (1-\varepsilon_j)
\longrightarrow 0$ and that $\sum_n \d_n < +\infty$, letting first $n
\ra +\infty$ and then $m \ra +\infty$ in \eqref{app-y_n} the
conclusion follows.

\qed
\\

\section{Some definitions and known results}

\subsection{Stable convergence and its variants}

 We recall here some basic definitions and results. For more details,
 we refer the reader to \cite{cri-let-pra-2007, hall-1980} and the
 references therein.\\

\indent Let $(\Omega, {\cal A}, P)$ be a probability space, and let
$S$ be a Polish space, endowed with its Borel $\sigma$-field. A {\em
  kernel} on $S$, or a random probability measure on $S$, is a
collection $K=\{K(\omega):\, \omega\in\Omega\}$ of probability
measures on the Borel $\sigma$-field of $S$ such that, for each
bounded Borel real function $f$ on $S$, the map
$$
\omega\mapsto
K\!f(\omega)=\int f (x)\, K(\omega)(dx)
$$
is $\cal A$-measurable. Given a sub-$\sigma$-field $\mathcal H$ of
$\cal A$, a kernel $K$ is said $\mathcal H$-measurable if all the
above random variables $K\!f$ are $\mathcal H$-measurable.\\

\indent On $(\Omega, {\mathcal A},P)$, let $(Y_n)$ be a sequence of
$S$-valued random variables, let $\mathcal H$ be a sub-$\sigma$-field
of $\cal A$, and let $K$ be a $\mathcal H$-measurable kernel on
$S$. Then we say that $Y_n$ converges {\em $\mathcal H$-stably} to
$K$, and we write $Y_n\stackrel{{\mathcal H}-stably}\longrightarrow
K$, if
$$
P(Y_n \in \cdot \,|\, H)\stackrel{weakly}\longrightarrow
E\left[K(\cdot)\,|\, H \right]
\qquad\hbox{for all } H\in{\mathcal H}\; \hbox{with } P(H) > 0.
$$
In the case when ${\mathcal H}={\mathcal A}$, we simply say that $Y_n$
converges {\em stably} to $K$ and we write
$Y_n\stackrel{stably}\longrightarrow K$. Clearly, if
$Y_n\stackrel{{\mathcal H}-stably}\longrightarrow K$, then $Y_n$
converges in distribution to the probability distribution
$E[K(\cdot)]$. Moreover, the $\mathcal H$-stable convergence of $Y_n$
to $K$ can be stated in terms of the following convergence of
conditional expectations:
\begin{equation}\label{def-stable}
E[f(Y_n)\,|\, {\cal H}]\stackrel{\sigma(L^1,\, L^{\infty})}\longrightarrow
K\!f
\end{equation}
for each bounded continuous real function $f$ on $S$. \\

\indent In \cite{cri-let-pra-2007} the notion of $\cal H$-stable
convergence is firstly generalized in a natural way replacing in
(\ref{def-stable}) the single sub-$\sigma$-field $\cal H$ by a
collection ${\cal G}=({\cal G}_n)$ (called conditioning system) of
sub-$\sigma$-fields of $\cal A$ and then it is strengthened by
substituting the convergence in $\sigma(L^1,L^{\infty})$ by the one in
probability (i.e. in $L^1$, since $f$ is bounded). Hence, according to
\cite{cri-let-pra-2007}, we say that $Y_n$ converges to $K$ {\em
  stably in the strong sense}, with respect to ${\cal G}=({\cal
  G}_n)$, if
\begin{equation}\label{def-stable-strong}
E\left[f(Y_n)\,|\,{\cal G}_n\right]\stackrel{P}\longrightarrow K\!f
\end{equation}
for each bounded continuous real function $f$ on $S$.\\

\indent Finally, a strengthening of the stable convergence in the
strong sense can be naturally obtained if in (\ref{def-stable-strong}) we
replace the convergence in probability by the almost sure convergence:
given a conditioning system ${\cal G}=({\cal G}_n)$, we say that $Y_n$
converges to $K$ in the sense of the {\em almost sure conditional
  convergence}, with respect to ${\cal G}$, if
\begin{equation}\label{def-as-cond}
E\left[f(Y_n)\,|\,{\cal G}_n\right]\stackrel{a.s.}\longrightarrow K\!f
\end{equation}
for each bounded continuous real function $f$ on $S$. Evidently, this
last type of convergence can be reformulated using the conditional
distributions. Indeed, if $K_n$ denotes a version of the conditional
distribution of $Y_n$ given ${\mathcal G}_n$, then the random variable
$K_n\!f$ is a version of the conditional expectation
$E\left[f(Y_n)|{\cal G}_n\right]$ and so we can say that $Y_n$ converges
to $K$ in the sense of the almost sure conditional convergence, with
respect to $\mathcal F$, if, for almost every $\omega$ in $\Omega$,
the probability measure $K_n(\omega)$ converges weakly to
$K(\omega)$. The almost sure conditional convergence has been
introduced in \cite{crimaldi-2009} and, subsequently, employed by
others in the urn model literature (e.g. \cite{aletti-2009, z}).
\\

We now conclude this section with a convergence result for martingale
difference arrays. \\ Given a conditioning system
$\mathcal{G}=(\mathcal{G}_n)_n$, if $\mathcal{U}$ is a
sub-$\sigma$-field of $\mathcal{A}$ such that, for each real
integrable random variable $Y$, the conditional expectation
${\mathrm{E}}[Y\,|\,{\mathcal{G}}_n]$ converges almost surely to the
conditional expectation ${\rm \mathrm{E}}[Y\,|\,\mathcal{U}]$, then we
shall briefly say that $\mathcal{U}$ is an {\em asymptotic
  $\sigma$-field} for $\mathcal{G}$. In order that there exists an
asymptotic $\sigma$-field $\mathcal{U}$ for a given conditioning
system $\mathcal{G}$, it is obviously sufficient that the sequence
$(\mathcal{G}_n)_n$ is increasing or decreasing.  (Indeed we can take
$\mathcal{U}=\bigvee_n\mathcal{G}_n$ in the first case and
${\mathcal{U}}=\bigcap_n\mathcal{G}_n$ in the second one.)

\begin{theorem}\label{fam_tri_vet_as_inf}
(Theorem A.1 in \cite{crimaldi-2009}) \\ On~$(\Omega,\mathcal{A},P)$,
  for each $n\geq 1$, let $(\mathcal{F}_{n,h})_{h\in\N}$ be a
  filtration and $(M_{n,h})_{h\in\N}$ a real martingale with respect
  to $({\mathcal{F}}_{n,h})_{h\in\N}$, with $M_{n,0}=0$, which
  converges in $L^1$ to a random variable $M_{n,\infty}$. Set

\begin{equation*}
X_{n,j}:=M_{n,j}-M_{n,j-1}\quad\hbox{for } j\geq 1,\quad
U_n:=\textstyle\sum_{j\geq 1}X_{n,j}^2,\quad
X_n^*:=\textstyle\sup_{j\geq 1}\;|X_{n,j}|.
\end{equation*}

Further, let $(k_n)_{n\geq 1}$ be a sequence of strictly positive
integers such that $k_nX_n^*\stackrel{a.s.}\to 0$ and let
$\mathcal{U}$ be a sub-$\sigma$-field which is asymptotic for the
conditioning system $\mathcal{G}$ defined by
$\mathcal{G}_n={\mathcal{F}}_{n,k_n}$. Assume that the sequence
$(X_n^*)_n$ is dominated in $L^1$ and that the sequence $(U_n)_n$
converges almost surely to a positive real random variable $U$ which
is measurable with respect to $\mathcal{U}$.  \\

Then, with respect to the conditioning
system $\mathcal{G}$, the sequence $(M_{n,\infty})_{n}$ converges to the
Gaussian kernel ${\mathcal{N}}(0,U)$ in the sense of the almost sure
conditional convergence.
\end{theorem}

\subsection{Durrett-Resnick result}

We recall the following convergence result for martingale difference arrays:

\begin{theorem}\label{dur-res-th}
(Th.~2.5 in \cite{dur-res-1978})\\ Let $(X_{n,k})$ be a
  square-integrable martingale difference array with respect to
  $(\mathcal{F}_{n,k})$. Suppose that $(\mathcal{F}_{n,k})$ increases as $n$
  increases and let $k_n(t)$ a non-decreasing right continuous
  function with values in $\N$ such that the following conditions hold
  true: \\[2pt] (a1) for each $t>0$,
$$
\sum_{k=1}^{k_n(t)}{\mathrm{E}}[X_{n,k}^2\,|\,\mathcal{F}_{n,k-1}]
\stackrel{P}\longrightarrow V_t
$$
where $P(t\mapsto V_t\; \hbox{ is continuous})=1$;
\\[5pt]
(b1) for each $\epsilon>0$,
$$
\sum_{k=1}^{k_n(1)}
{\mathrm{E}}[X_{n,k}^2I_{\{|X_{n,k}|>\epsilon \}}\,|\,\mathcal{F}_{n,k-1}]
\stackrel{P}\longrightarrow 0.
$$
Then, if we set $S_{n,k_n(t)}=\sum_{k=1}^{k_n(t)}X_{n,k}$, we have
$$
S^{(n)}=(S_{n,k_n(t)})_{t\geq 0}
\stackrel{d}\longrightarrow
\widetilde{W}=\big(W_{V_t}\big)_{t\geq 0}
\qquad\hbox{(w.r.t. Skorohod's topology)},
$$
where $W=(W_t)_{t\geq 0}$ is a Wiener process
independent of $V=(V_t)_{t\geq 0}$.
\end{theorem}

\begin{remark}\label{remark-cond}
\rm
If there exists a number $u\geq 2$ such that
\be{suff}
\sum_{k=1}^{k_n(1)}
{\mathrm{E}}[\,|X_{n,k}|^u\,]\longrightarrow 0,
\ee
then condition (b1) holds true with convergence in $L^1$.
Indeed, it is enough to observe that
\begin{equation*}
\sum_{k=1}^{k_n(1)}
{\mathrm{E}}[X_{n,k}^2I_{\{|X_{n,k}|>\epsilon\}}]=
\sum_{k=1}^{k_n(1)}
{\mathrm{E}}[|X_{n,k}|^u |X_{n,k}|^{-(u-2)}I_{\{|X_{n,k}|>\epsilon\}}]
\leq
\frac{1}{\epsilon^{u-2}}
\sum_{k=1}^{k_n(1)}{\mathrm{E}}[\,|X_{n,k}|^u\,]
\longrightarrow 0.
\end{equation*}
\end{remark}

\subsection{H\'ajek-R\'enyi inequality}

We recall the following martingale inequality (e.g. \cite {cho-tei}):
\\

\begin{theorem}
If $M_n=\sum_{j=1}^n\xi_j$ is a square integrable martingale and
$(a_n)$ is a positive, nondecreasing sequence of numbers, then for
each $\lambda>0$, we have
$$
P\left(\max_{1\leq j\leq n}\frac{|M_j|}{a_j}\geq \lambda\right)
\leq
\frac{1}{\lambda^2}\sum_{j=1}^n\frac{\mathrm{E}[\xi_j^2]}{a_j^2}\,.
$$
\end{theorem}

\subsection{Barbour's transform}

Let $D=D[0,+\infty)$ be the space of right-continuous functions with
left limits on $[0,+\infty)$, endowed with  the classical Skorohod's
topology (e.g. \cite{Bill2}).
\\[5pt]\noindent
Let $T$ be the subspace of functions $f(t)$ in the space $D$
such that
\begin{equation}\label{condA}
\limsup_{t\to +\infty}\frac{|f(t)|}{t}=0
\end{equation}
\begin{equation}\label{condB}
\int_1^{+\infty}\frac{|f(t)|}{t^2}\,{\rm d}t<+\infty
\end{equation}
\begin{equation}\label{condC}
\int_0^1\frac{|f(t)|}{t}\,{\rm d}t<+\infty.
\end{equation}
Let $m$ be the metric on $T$ such that $m(f_1,f_2)$ is the infimum of
those $\epsilon>0$ for which there exists some continuous strictly
increasing function $\lambda:[0,+\infty)\mapsto[0,+\infty)$ with
$\lambda(0)=0$, such that
\begin{equation}\label{condAA}
\sup_{t\geq 0}
\frac{f_1(t)-f_2(\lambda(t))}{t+1}<\epsilon
\end{equation}
\begin{equation}\label{condBB}
\int_1^{+\infty}\frac{|f_1(t)-f_2(\lambda(t))|}{t^2}\,{\rm d}t<\epsilon
\end{equation}
\begin{equation}\label{condCC}
\int_0^1\frac{|f_1(t)-f_2(\lambda(t))|}{t}\,{\rm d}t<\epsilon
\end{equation}
\begin{equation}\label{condDD}
\sup_{t \neq s} \left| \log \frac{\lambda(t) - \lambda(s)}{t-s} \right| <
\epsilon.
\end{equation}

Let $T_1$ and $m_1$ be defined similarly, without the restrictions
(\ref{condC}) and (\ref{condCC}). We shall denote by $T^*$ and $T_1^*$
the corresponding subspaces of the space $D^*=D^*[0,+\infty)$ of
  left-continuous functions with right limits on $[0,+\infty)$
    (endowed with the corresponding Skorohod's topology).  \\[5pt]

The topology induced by $m$ on $T$ is stronger than the Skorohod's
topology. Moreover, the {\em Barbour's transform} $g:T\to T_1^*$
defined as
$$
g(f)(0)=0\quad\hbox{and}\quad
g(f)(t):=\int_{1/t}^{+\infty} s^{-1}{\rm d}f(s)=
-tf(t^{-1})+\int_{1/t}^{+\infty} s^{-2}f(s)\,{\rm d}s
\qquad\hbox{for } t\in(0,+\infty).
$$

is continuous and, if $W$ is a Wiener process, then also $g(W)$ is
a Wiener process. Finally, the following result holds (for more
details, see \cite{gouet, hey, mul, whi}).

\begin{theorem}\label{gouet-th}
Let $(Y^{(n)})_n$ be a sequence of stochastic processes satisfying the
following conditions:
\begin{itemize}
\item[(a2)] $Y^{(n)}\stackrel{d}\longrightarrow \widetilde{W}$
  (w.r.t. Skorohod's topology), where $\widetilde W$ is a stochastic
  process of the form $\widetilde{W}_t=W_{V_t}$ where $W$ is a Wiener
  process and $V$ is a stochastic process, independent of $W$ and with
  $P(t\mapsto V_t\; \hbox{ is continuous})=1$;
\item[(b2)] for each $n$ and $\epsilon>0$,
  $Y_n(t)=o(t^{1/2+\epsilon})$ a.s. as $t\to +\infty$;
\item[(c2)] for each $\theta>1/2$, $\epsilon>0$ and
$\eta>0$, there exists $t_0$ such that
\begin{equation}\label{inequality-th}
P\left\{
\sup_{t\geq t_0}\frac{ |Y^{(n)}_t| }{ t^{\theta}}
>\epsilon
\right\}
\leq\eta.
\end{equation}
\end{itemize}
Then each $Y^{(n)}$ takes values in the space $T$ and
$Y^{(n)}\stackrel{d}\longrightarrow \widetilde{W}$ on $(T,m)$.
\end{theorem}


\begin{thebibliography}{99}
\bibitem{aletti-2009} Aletti G., May C. and Secchi P. (2009).
A central limit theorem, and related results, for a two-color randomly
reinforced urn. Adv. Appl. Prob. 41, 829-844.
\bibitem{aletti} Aletti, G. and Ghiglietti, A. (2016). Interacting
  Generalized P\'olya Urn Systems. arXiv preprint arXiv:1601.01550.
\bibitem{bai-et-al} Bai Z.~D., Hu F. and Zhang L.~X. (2002). Gaussian
  approximation theorems for urn models and their
  applications. Ann. Appl. Probab., {\bf 12}(4), 1149-1173.
  \bibitem{bar} Barab\'asi, A. L., and Albert, R. (1999). Emergence of
  scaling in random networks. Science, 286(5439), 509-512.
\bibitem{bas-das} Basak G.~K. and Dasgupta A. (2005). A functional
  central limit theorem for a class of urn models. Proc. Indian
  Acad. Sci. (Math. Sci.), {\bf 115}(4), 493-498.
\bibitem{Bena} Bena\"im, M. (1999). Dynamics of stochastic
  approximation algorithms. In J. Az\'ema, M. \'Emery, M. Ledoux,
  M. Yor (Eds.), 1-68.  Springer Berlin Heidelberg.
   \bibitem{ben} Bena\"im, M., Benjamini, I., Chen, J., and Lima,
   Y. (2015). A generalized P\'{o}lya's urn with graph based
   interactions. Random Structures \& Algorithms, 46(4), 614-634.
\bibitem{Bill2} Billingsley P. (1999) {Convergence of Probability
  Measures (second edition)}.  Wiley, New York.
\bibitem{cho-tei} Chow Y.~S. and Teicher H. (1988). {Probability
  theory: independence, interchangeability, martingales}. Springer,
  New York.
  \bibitem{collet} Collet, F. (2014). Macroscopic limit of a bipartite
    Curie-Weiss model: A dynamical approach. Journal of Statistical
    Physics, 157(6), 1301-1319.
\bibitem{daipra15} Collet F., Dai Pra P. and Formentin
   M. (2015). Collective periodicity in mean-field models of
   cooperative behavior. Nonlinear Differential Equations and
   Applications NoDEA, 22(5), 1461-1482.
   \bibitem{contucci} Contucci, P., and Ghirlanda, S. (2007). Modeling
  society with statistical mechanics: an application to cultural
  contact and immigration. Quality \& quantity, 41(4), 569-578.
  \bibitem{cri-dai-min} Crimaldi I., Dai Pra P., Minelli I.G. (2016),
  Fluctuation theorems for synchronization of interacting P\'olya's
  urns, Stoch. Proc. Appl., {\bf 126} (3), 930-947.
\bibitem{crimaldi-2009} Crimaldi I. (2009).
{An almost sure conditional convergence result and an application to a
generalized P\'olya urn}. Internat. Math. Forum, {\bf 4}(23), 1139-1156.
\bibitem{cri-let-pra-2007} Crimaldi I., Letta G. and Pratelli L. (2007).
A strong form of stable convergence. In S\'eminaire de
Probabilit\'es XL (Lecture Notes Math. 1899), Springer, Berlin, 203-225.
\bibitem{daipra-2014} Dai Pra P., Louis P.-Y. and Minelli I.~G. (2014).
  Syncronization via interacting reinforcement. J. Appl. Probab.,
  {\bf 51}(2), 556-568.
\bibitem{davis} Davis, B. (1990). Reinforced random
     walk. Prob. Th. Rel. Fields, {\bf 84}(2), 203-229.
\bibitem{dur-res-1978} Durrett R. and Resnick S.~I. (1978).
{Functional limit theorems for
  dependent variables}. Ann. Probab., {\bf 6}(5), 829-846.
\bibitem{fabian} Fabian V. (1968). On asymptotic normality in
  stochastic approximation. Ann. Math. Statist., {\bf 39}(4),
  1327-1332.
  \bibitem{fagnani} Como G. and Fagnani F. (2011). Scaling limits
    for continuous opinion dynamics systems. The Annals of Applied
    Probability, 21(4), 1537-1567.
\bibitem{fr65} Freedman D. A. (1965). Bernard Friedman's urn. The
  Annals of Mathematical Statistics, 956-970.
\bibitem{GB} Giacomin G., Pakdaman K., Pellegrin X., and Poquet C
  (2012).  Transitions in Active Rotator Systems: Invariant Hyperbolic
  Manifold Approach. SIAM J. Math. Anal., 44(6), 4165-4194.
 \bibitem{gibbs-2002} Gibbs A.~L. and Su F.~E. (2002). On choosing and
  bounding probability metrics.  International Statistical Review,
  {\bf 70}(3), 419-435.
\bibitem{gouet} Gouet R. (1993).
{Martingale functional central limit theorems for a
  generalized P\'olya urn}. Ann. Probab., {\bf 21}(3), 1624-1639.
\bibitem{hall-1980} Hall P. and Heyde C.~C. (1980).
Martingale Limit Theory and Its Applications. Academic Press, New York.
\bibitem{hey} Heyde C.~C. (1977).
{On central limit and iterated logarithm supplements
to the martingale convergence theorem}. J. Appl. Prob.,
{\bf 14}, 758-775.
\bibitem{hida} Hida T. (1980) Brownian motion. Springer, New York.
 \bibitem{remco} Hofstad (van der), R. (2009). Random graphs and
    complex networks. Available on
    http://www. win. tue. nl/rhofstad/NotesRGCN.
\bibitem{janson} Janson S. (2004). Functional limit theorems for
  multitype branching processes and generalized P\'olya
  urns. Stoch. Proc. Appl., {\bf 110}, 177-245.
\bibitem{Kush} Kushner H. and Yin G. (1997). Stochastic Approximation
  Algorithms and Applications, Springer.
\bibitem{launay-2012} Launay M. and Limic V. (2012).
Generalized interacting urn models. arXiv
preprint, arXiv:1207.5635.
\bibitem{launay} Launay M. (2011). Interacting urn models.
arXiv preprint, arXiv:1101.1410.
\bibitem{mah} Mahmoud, H. (2008). P\'olya urn models. CRC press.
\bibitem{marsili} Marsili M. and Valleriani
  A. (1998). Self-organization of interacting P\'olya urns.
Eur. Phys. J. B 3 (4), 417-420.
\bibitem{mul} Muller D.~W. (1968).  {Verteilungs-Invarianzprinzipien
  f\"ur das starke Gesetz der grossen Zahl.}
  Z. Wahrsch. Verw. Gebiete, {\bf 10}, 173-192.
\bibitem{paganoni} Paganoni A. M. and Secchi P. (2004).  Interacting
 reinforced-urn systems. Adv. Appl. Probab., 36(3), 791-804.
\bibitem{perth1} Pakdaman K., Perthame B. and  Salort D. (2010).
Dynamics of a structured neuron population.
Nonlinearity, 23(1), 55-75.
\bibitem{pemantle-2007} Pemantle R. (2007).
A survey of random processes with reinforcement. Prob. Surveys 4, 1-79.
\bibitem{rob-sie} Robbins H. and Siegmund D. (1971). A convergence
  theorem for non negative almost supermartingales and some
  applications. Optimizing Methods in Statistics. Chapter
  Herbert Robbins Selected Papers, 111-135, Springer.
  \bibitem{neeraja} Sahasrabudhe N.  (2015). Synchronization and
    Fluctuation Theorems for Interacting Friedman Urns. Preprint.
 \bibitem{whi} Whitt W. (1972). {Stochastic abelian and tauberian theorems}.
Z. Wahrsch. Verw. Gebiete, {\bf 22}, 251-267.
\bibitem{z} Zhang L. X. (2014). A Gaussian process approximation for
  two-color randomly reinforced urns. Electron. J. Probab., {\bf 19}(86),
  1-19.
\end{thebibliography}
\end{document}